\documentclass{article}

\usepackage{proof}
\usepackage{latexsym}
\usepackage{amsfonts}
\usepackage{amssymb}
\usepackage{cite}
\usepackage{color}

\newcommand{\blem}{\begin{lemma}}
\newcommand{\elem}{\end{lemma}}
\newcommand{\bth}{\begin{theorem}}
\newcommand{\ethm}{\end{theorem}}
\newcommand{\benu}{\begin{enumerate}}
\newcommand{\eenu}{\end{enumerate}}
\newcommand{\bdes}{\begin{description}}
\newcommand{\edes}{\end{description}}
\newcommand{\bdf}{\begin{definition}}
\newcommand{\edf}{\end{definition}}
\newcommand{\bcor}{\begin{cor}}
\newcommand{\ecor}{\end{cor}}
\newcommand{\bprp}{\begin{proposition}}
\newcommand{\eprp}{\end{proposition}}
\newcommand{\bmlem}{\begin{mlemma}}
\newcommand{\emlem}{\end{mlemma}}
\newcommand{\bclm}{\begin{claim}}
\newcommand{\eclm}{\end{claim}}
\newcommand{\bprf}{{\bf Proof}.\hspace{2mm}}

\newcommand{\eprf}{\hspace*{\fill} $\Box$}

\newcommand{\ovl}{\overline}
\newcommand{\beqn}{\begin{equation}}
\newcommand{\eeqn}{\end{equation}}
\newcommand{\beqnarr}{\begin{eqnarray}}
\newcommand{\eeqnarr}{\end{eqnarray}}
\newcommand{\beqnarrs}{\begin{eqnarray*}}
\newcommand{\eeqnarrs}{\end{eqnarray*}}
\newcommand{\ull}{\underline{\ll}}

\newcommand{\spand}{\,\&\,}

\newcommand{\Natural}{\mathbb{N}}

\newtheorem{theorem}{Theorem}[section]
\newtheorem{definition}[theorem]{Definition}
\newtheorem{proposition}[theorem]{Proposition}
\newtheorem{lemma}[theorem]{Lemma}
\newtheorem{cor}[theorem]{Corollary}
\newtheorem{mlemma}[theorem]{Main Lemma}
\newtheorem{claim}[theorem]{Claim}

\newcommand{\alp}{\alpha}
\newcommand{\eps}{\epsilon}

\newcommand{\del}{\delta}
\newcommand{\Del}{\Delta}
\newcommand{\ome}{\omega}
\newcommand{\Ome}{\Omega}
\newcommand{\bet}{\beta}
\newcommand{\gam}{\gamma}
\newcommand{\Gam}{\Gamma}
\newcommand{\kap}{\kappa}
\newcommand{\sig}{\sigma}
\newcommand{\Sig}{\Sigma}
\newcommand{\tht}{\theta}

\newcommand{\Lam}{\Lambda}
\newcommand{\vphi}{\varphi}

\newcommand{\fal}{\forall}
\newcommand{\exi}{\exists}
\newcommand{\rarw }{\rightarrow}
\newcommand{\Rarw }{\Rightarrow}

\newcommand{\lrarw}{\leftrightarrow}
\newcommand{\Lrarw}{\Leftrightarrow}
\newcommand{\uarw}{\uparrow}
\newcommand{\darw}{\downarrow}

\newcommand{\cala}{{\cal A}}
\newcommand{\calb}{{\cal B}}

\newcommand{\calg}{{\cal G}}
\newcommand{\calh}{{\cal H}}

\newcommand{\calt}{{\cal T}}

\newcommand{\calL}{{\cal L}}

\newcommand{\calP}{{\cal P}}

\newcommand{\rk}{\mbox{{\rm rk}}}

\newcommand{\msfiv}{\mbox{\hspace{5mm}}}

\newcommand{\setm}{\setminus}

\title{Hydra games for recursively Mahlo operations}
\author{Toshiyasu Arai
\\
\\
Graduate School of Science,
Chiba University
\\
1-33, Yayoi-cho, Inage-ku,
Chiba, 263-8522, JAPAN
\\
tosarai@faculty.chiba-u.jp
}
\date{}
\begin{document}
\maketitle
\begin{abstract}
Encouraged by W.  Buchholz \cite{Buchholz98}, 
a hydra game is proposed,  and 
the fact that every hydra eventually die out is shown to be equivalent (over a weak arithmetic)
 to the 1-consistency of set theory
{\sf KPM} for recursively Mahlo universes.
\end{abstract}


\section{Introduction}\label{sect:intro}
In M. Rathjen\cite{Rathjen91}, W. Buchholz\cite{Buchholz90} and \cite{odMahlo,ptMahlo} 
the set theory {\sf KPM} for recursively Mahlo universes has been analyzed proof-theoretically.

As to the proof-theoretic analyses on such strong impredicative theories, 
let us quote from Buchholz\cite{Buchholz98}:
\begin{quote}
Contemporary ordinal-theoretic proof  theory (i.e., the part of proof theory concerned 
with ordinal analyses of strong impredicative theories) suffers from the extreme 
(and as it seems unavoidable)
complexity and opacity of its main tool, the  ordinal notation systems.
This is not only a technical stumbling block which prevents most proof-theorists 
from a closer engagement in that field, but it also calls the achieved results into question, 
at least as long as these results do not have interesting consequences, such as e.g., 
foundational reductions or intuitively graspable combinatorial independence results.
\end{quote}

 If proofs or constructions looks too complicated
\footnote{Indeed, it's complicated as  compared with those for predicative theories such as {\sf PA}.}
to grasp,
and this makes us doubtful about what  we have gained,
I would reply that this defect are mainly due to the scarcity of our experiences of mathematics in 
the strength of strong impredicative theories $T$.

One thing we can do is to give alternative proofs, thereby could shed light on the 
same results from another angle, and gain an insight in mathematical reasoning and structures codified in $T$.
Another thing to be done is to find combinatorial independence results.
This line of research was suggested and encouraged by W.  Buchholz \cite{Buchholz98}.
In a sense, such an (optimal) independence result might be viewed as a finitary essence of $T$.
One prototype of combinatorial independence results is the hydra games in 
Kirby-Paris\cite{KS} and Buchholz\cite{Buchholz87}.
This says that given a hydra game, a theory such as {\sf PA} or $(\Pi^{1}_{1}{\rm -CA})+{\rm BI}$ proves that
\textit{each}  hydra eventually die out, but the theory in question does not prove its universal  closure,
\textit{any} hydra must die out.
\\

In this paper
a hydra game for recursively Mahlo ordinals is proposed,  
 and a result of the same kind
is shown for the games, and {\sf KPM}.


A tree $(T,<)$ is said to be  \textit{structured} if
the (finite) set of immediate successors $\{s: t<s \spand \lnot\exi u(t<u<s)\}$ of each node $t$ in $T$
is linearly ordered.

A \textit{hydra} is a triple $(T,<,\ell)$ such that $(T,<)$ is a finite and structured tree
and $\ell: T\to \{\star\}\cup Lb_{0}\cup Lb_{1}$,
where $\star$ is the label attached to the root of the tree $T$,
$Lb_{0}=\{1,n\cdot m,n\cdot A,n\cdot *_{\ome}, n\cdot*_{\mu}: 0<n,m<\ome, A\in L\}$
is the set of labels for leaves, and
$Lb_{1}=\{\ome\}\cup\{\{A\}: A\in L\cup\{*_{\mu},\mu\}\}\cup
\{\vphi_{\mbox{{\scriptsize \boldmath$A$}}+n},D_{\mbox{\scriptsize \boldmath$A$}}:\mbox{\boldmath$A$}\in L^{*},n<\ome\}$
is the set of labels for internal nodes.
 $L$ is the set of \textit{labels} defined below, and
 $L^{*}$ the set of finite sets of labels in $L$.
 Each hydra and each label in $L$ is a term over symbols $\{n: 0<n<\ome\}\cup\{\cdot,\ome,\mu, D, \{\,\}, *_{\ome},*_{\mu},\vphi,d\}$.
The set $L$ of labels ordered by a linear order $<$ with the largest element $\mu$.
The set of hydras $\calh$ and the set of labels $L$ are defined simultaneously.

Hydras produce a finite set of labels, as the game goes.
In some limit cases of the hydra game, 
a hydra $(T,<,\ell)$ freely chooses a label from the finite set of labels, which are available for the current hydra.
The set of their labels might grow in some cases called (Production).

A free choice of labels means that,
for a hydra $H$, a finite set $lb$ of labels and natural numbers $\ell$,
there are finitely many \textit{possible moves} written as $(H,lb)\to_{\ell}(K,lb^{\prime})$.
Given a hydra $H_{0}$ and a finite set $lb_{0}$ of labels in $L$, 
a finitely branching tree is obtained as follows.
For $t\in{}^{<\ome}\ome$
we define moves $(H_{0}[t], lb[t])$ in the hydra game.
$H_{0}[\eps]:=H_{0}$ and $lb[\eps]:=lb_{0}$ for the empty sequence $\eps$.
$\{(H_{0}[t*(i)], lb[t*(i)])\}_{i}$ is the set of the pairs $(K,lb^{\prime})$ such that
$(H_{0}[t],lb[t])\to_{\ell}(K,lb^{\prime})$ for the length $\ell=|t|$ of the finite sequence $t$.
The finitely branching tree 
$Tr(H_{0},lb_{0})=\{t\in {}^{<\ome}\ome: (H_{0}[t],lb_{0}[t]) \mbox{ {\rm is defied}}\}$
is thus obtained from $H_{0}$ and $lb_{0}$.
The tree $Tr(H_{0},lb_{0})$ is 
seen to be well founded for every hydra $H_{0}$ and every finite set $lb_{0}$ of labels.
Let $F[H_{0},lb_{0}]$ denote the length of maximal runs in the game, i.e., the height of the tree:
\[
F[H_{0},lb_{0}]=\max\{|t|: t\in{}^{<\ome}\ome,  H_{0}[t] \mbox{ {\rm is defined with }}
lb[\eps]=lb_{0}
\}
.\]
Let $H+n:=H+1+\cdots+1$ be the hydra obtained from the hydra $H$ by adding
the trivial hydra $1$ in $n$-times.
Now our theorem runs as follows.

\bth\label{th:main}

\benu
\item
\label{th:main1}
Each provably total recursive function in {\sf KPM} is dominated by a function
$n\mapsto F[H_{0}+n,\emptyset]$ for some hydra $H_{0}$.
\item
\label{th:main2}
Conversely for each hydra $H_{0}$ and each finite set $lb_{0}\subset L$, 
the function $n\mapsto F[H_{0}+n,lb_{0}]$
is provably total recursive function in {\sf KPM}.
\item
\label{th:main3}
The fact that for every 
hydra $H_{0}$, the hydra game eventually terminates, i.e.,
 the tree $\{t\in{}^{<\ome}\ome:  H_{0}[t] \mbox{ {\rm is defined with }} lb[\eps]=\emptyset\}$ is finite,
 or equivalently the $\Pi^{0}_{2}$-statement
 $\fal H_{0}\, F[H_{0},\emptyset]\darw$
  is equivalent to the 1-consistency $\mbox{{\rm RFN}}_{\Pi^{0}_{2}}(\mbox{{\sf KPM}})$ of {\sf KPM}
  over the
 elementary arithmetic {\sf EA}.
\eenu
\end{theorem}

Let us mention the contents of the paper.
In section \ref{sect:hydragame} the hydra game is defined
through a linear ordering $A<B$ on labels, which is based on
an assignment of ordinal diagrams $o(H),o(A)\in O(\mu)$ to hydras $H$ and labels $A\in L$.
In section \ref{sect:prov} we show that
$d_{\Ome}\left(o(K)\#o(lb^{\prime})\right)< d_{\Ome}\left(o(H)\#o(lb)\right)$
when $(H,lb)\to_{\ell}(K,lb^{\prime})$ is a possible move.
Thus Theorem \ref{th:main}.\ref{th:main2} follows from the fact in \cite{odMahlo} that
the wellfoundedness up to each ordinal diagram$<\Ome$ is provable in {\sf KPM}.
In section \ref{sect:unprov} we introduce first
a theory $[\Pi^{0}_{1},\Pi^{0}_{1}]$-Fix for non-monotonic inductive definitions
of $[\Pi^{0}_{1},\Pi^{0}_{1}]$-operators in \cite{Richter-Aczel74}.
In \cite{wienpi3} it is shown that the 1-consistency of the set theory {\sf KPM} is reduced to one of
the theory $[\Pi^{0}_{1},\Pi^{0}_{1}]$-Fix.
Second we assign hydras to proofs in the theory.
In section \ref{sec:rewriting}
we define rewritings on proofs in such a way that each rewriting corresponds to 
a move on hydras attached to proofs.
Theorems \ref{th:main}.\ref{th:main1} and \ref{th:main}.\ref{th:main3} are concluded.
Finally the linearity of the relation $A<B$ on labels is briefly discussed.

\section{Hydra game}\label{sect:hydragame}

In this section we introduce a hydra game for recursively Mahlo ordinals.


In the next Definition \ref{df:hydra} the set $\calh=\calh_{0}\cup\calh_{1}$ of \textit{hydras} 
and the set $L$ of 
\textit{labels} are defined simultaneously.
Also we define a subset $\calt_{i}\subset\calh_{i}$ for $i=0,1$.

\bdf\label{df:hydra}{\rm (Hydras)}
 \benu
  \item\label{df:hydra.7}
 $L=\{d_{\mu}(h_{0}) :h_{0}\in\calh_{0}\}\cup\{d_{d_{\mu}(h_{0})}(h_{1}): h_{0}\in\calh_{0}, h_{1}\in\calh_{1}\}$.
 
 {\rm 
 $L^{*}$ denotes the set of all finite sets of labels in $L$.
 The singleton $\{A\}$ is identified with labels $A\in L$.
 $R(A)$ holds iff either $A=\mu$ or $A=d_{\mu}(h)$.}
 

  \item\label{df:hydra.-1}
  $\calt_{i}\subset\calh_{i}$ {\rm for} $i=0,1$.
  
 \item\label{df:hydra.0}
 {\rm $0\in\calh_{0}\cap\calh_{1}$
 and $1\in\calt_{0}\cap\calt_{1}$.
 }

  \item\label{df:hydra.1}
{\rm If  $h_{0},\ldots,h_{k}\in \calt_{i}\, (k>0)$, then $(h_{0}+\cdots +h_{k})\in \calh_{i}$ for $i=0,1$.
}

 \item\label{df:hydra.5}
{\rm $\{n\cdot *_{\ome}, n\cdot *_{\mu}\}\cup\{n\cdot m, n\cdot A: 0<n,m<\ome, A\in L\}\subset\calt_{0}$.


}

   \item\label{df:hydra.3}
     {\rm If $h\in \calh_{0}$,
     then $\ome(h)\in \calt_{0}$.
  }
 
  \item\label{df:hydra.6}
 {\rm If $h\in\calh_{0}$, then 
 $\{\mu\}(h),\{*_{\mu}\}(h)\in\calt_{0}$.
 Also if $h\in\calh_{i}$ and $A\in L$, then $\{A\}(h)\in\calt_{i}$ for $i=0,1$.
 }
 
\item\label{df:hydra.4}
{\rm 
If $h\in \calh_{0}$ and $\mbox{\boldmath$C$}\in L^{*}$, then $D(\mbox{\boldmath$C$};h)\in \calt_{1}$.
Let $D_{\mbox{\scriptsize \boldmath$C$}}(h):=D(\mbox{\boldmath$C$};h)$.
}

 \item\label{df:hydra.2}
  {\rm If $h\in \calh_{1}$,
  $\emptyset\neq\mbox{\boldmath$C$}\in L^{*}$ and $n<\ome$, 
  then $\vphi(\mbox{\boldmath$C$}+n;h)\in \calt_{1}$.
  
  Let
  $\vphi(\emptyset;h):=h$ when $\mbox{\boldmath$C$}=\emptyset$ and $n=0$, and
  $\vphi_{\mbox{\scriptsize \boldmath$C$}+n}(h):=\vphi(\mbox{\boldmath$C$}+n;h)$.
  }


  
    

\eenu

\edf

\bdf\label{df:Lb}
{\rm
The set of labels $Lb(h)$ and the \textit{fixed part} $(h)_{f}\subset Lb(h)$ 
for hydras $h\in\calh$ are defined recursively as follows.
\benu
\item
$Lb(0)=Lb(1)=Lb(n\cdot *_{\ome})=Lb(n\cdot *_{\mu})=Lb(n\cdot m)=\emptyset$, and
$(0)_{f}=(1)_{f}=(n\cdot *_{\ome})_{f}=(n\cdot *_{\mu})_{f}=(n\cdot m)_{f}=\emptyset$.
\item
$Lb(h_{0}+\cdots+h_{k})=\bigcup_{i\leq k}Lb(h_{i})$, and $(h_{0}+\cdots+h_{k})_{f}=\bigcup_{i\leq k}(h_{i})_{f}$.
\item
$Lb(n\cdot C)=\{C\}$ and $(n\cdot C)_{f}=\emptyset$ for $C\in L$.
\item
$Lb(D(\mbox{\boldmath$C$};h))=Lb(\vphi_{\mbox{\scriptsize \boldmath$C$}+n}(h))=\mbox{\boldmath$C$}\cup Lb(h)$ and 
$(D(\mbox{\boldmath$C$};h))_{f}=(\vphi_{\mbox{\scriptsize \boldmath$C$}+n}(h))_{f}=\mbox{\boldmath$C$}\cup (h)_{f}$
for $\mbox{\boldmath$C$}\in L^{*}$.
\item
$Lb(\{A\}(h))=\{A\}\cup Lb(h)$ and $(\{A\}(h))_{f}=(h)_{f}$ for $A\in L$ and $n<\ome$.
\item
$Lb(\ome(h))=Lb(\{\mu\}(h))=Lb(\{*_{\mu}\}(h))=Lb(h)$ and
$(\ome(h))_{f}=(\{\mu\}(h))_{f}=(\{*_{\mu}\}(h))_{f}=(h)_{f}$
\eenu
}
\edf


In \cite{odMahlo} a system $(O(\mu),<)$ of ordinal diagram, a computable system of ordinal notations is defined,
and it is shown that {\sf KPM} proves the wellfoundedness up to each $\alp<\Ome$.
Let us recall a slightly modified system $(O(\mu),<)$ briefly.
The set $O(\mu)$ is generated from $0$ and $\mu$ by the addition $+$, 
the fixed point free binary Veblen function $\vphi\alp\bet\, (\alp,\bet<\mu)$,
the exponential above $\mu$, $\ome^{\alp}\,(\alp>\mu)$,
and the collapsing function $d:(\sig,\alp)\mapsto d_{\sig}\alp$ for the regular diagram $\sig$, i.e.,
either $\sig=\mu$ or $\sig=d_{\mu}\bet$ for a $\bet$.
$R$ denotes the set of all regular diagrams, and $\Ome:=d_{\mu}0$.
$\sig,\tau,\kap,\rho,\ldots$ denote regular diagrams.
Each $d_{\sig}\alp$ is a strongly critical number.
Crucial definitions are as follows.
$\alp\prec\bet$ iff either $\alp=d_{\bet}\gam$ or $\alp=d_{d_{\bet}\gam}\del$ 
for some $\gam,\del$.
 $\alp\preceq\bet:\Lrarw(\alp\prec\bet\lor\alp=\bet)$.
The set $K_{\sig}\alp$ of subdiagrams of $\alp$ is defined as follows.

 \begin{enumerate}
\item 
$K_{\sig}0=K_{\sig}\mu=\emptyset$,
$K_{\sig} (\alpha_1+\cdots +\alpha_n)=\bigcup \{K_{\sig} \alpha_i:1\leq i\leq n \}$, {\rm and}
 $K_{\sig} \vphi\alpha\bet=K_{\sig} \alpha\cup K_{\sig} \bet$.
 \item $K_{\sig}d_\tau\alp=\left\{
        \begin{array}{ll}
        \{d_{\tau}\alp\} &  \tau\preceq\sig
        \\
         K_\sig\tau\cup K_\sig\alp & \sig<\tau
          \\
         K_\sig\tau  & \tau<\sig \, \& \, \tau\not\preceq\sig
        \end{array}
        \right.$
 \end{enumerate}

For $\sig\neq\tau$, $d_{\sig}\alp<d_{\tau}\bet$
 iff one of the following conditions holds:
\begin{enumerate}
    \item $\sig<\tau \spand (\sig\leq d_{\tau}\bet \mbox{ {\rm or} } d_{\sig}\alp\leq K_\sig d_{\tau}\bet)$. 
    \item $\tau<\sig \spand d_{\sig}\alp<\tau \spand K_\tau d_{\sig}\alp<d_{\tau}\bet$.
 \end{enumerate}
$d_{\sig}\alp<d_{\sig}\bet$ iff one of the following conditions holds:
 \benu
 \item $d_{\sig}\alp\leq K_{\sig}\bet$.
 \item $K_{\sig}\alp<\bet \spand d_{\tau}\alp<d_{\tau}\bet$.
 \eenu
where 
$\tau=\min\{\tau\in R\cup\{\infty\}: (\sig<\tau<\infty \& K_\tau\{\alp,\bet\}\neq\emptyset) \mbox{ {\rm or} } \tau=\infty\}$,
and, by definition, $d_\infty\alp:=\alp$, $\fal\alp\in O(\mu)(\alp<\infty)$
and $\infty\not\in O(\mu)$.

We associate an ordinal diagram $o(h)\in O(\mu)$ for hydras $h$.

\bdf
{\rm 
We associate $o(A), o(h)\in O(\mu)$ for labels $A\in L$ and hydras $h\in\calh$ as follows.
\benu
\item
$o(d_{\mu}(h))=d_{\mu}(o(h))$ and 
$o(d_{d_{\mu}(h)}(h_{1}))=d_{o(d_{\mu}(h))}(o(d_{\mu}(h))\# o(h_{1}))$
with the natural (commutative) sum $\#$ in $O(\mu)$.

\item
$o(0)=0$ and $o(1)=1:=\ome^{0}$.

\item
$o(h_{0}+\cdots+h_{k})=o(h_{0})\#\cdots \# o(h_{k})$.

\item
$o(n\cdot *_{\ome})=\ome$, $o(n\cdot *_{\mu})=\mu$, $o(n\cdot m)=n\cdot m$ and $o(n\cdot A)=o(A)$.

\item
$o(\vphi_{\mbox{\scriptsize \boldmath$C$}+n}(h))=\vphi_{o(\mbox{\scriptsize \boldmath$C$})+n+1}(o(h))$,
where $o(\{C_{1},\ldots,C_{n}\})=o(C_{1})\#\cdots\#o(C_{n})$.

\item 
$o(\ome(h))=\ome^{o(h)}$.

\item
$o(\{\mu\}(h))=o(\{*_{\mu}\}(h))=\mu\# o(h)$.

\item
$o(\{A\}(h))=
\left\{
\begin{array}{ll}
o(A)\#1\#o(h) & \mbox{ if } h\in\calh_{0}
\\
\vphi_{o(A)}(o(h)) & \mbox{ if } h\in\calh_{1}
\end{array}
\right.
$.

\item
$o(D(\mbox{\boldmath$C$};h))=d_{\mu}(o(\mbox{\boldmath$C$})\# o(h))$.
\eenu

For $A,B\in L\cup\{0,\mu\}\cup\calh$, $n,m<\ome$, and 
$\mbox{\boldmath$A$},\mbox{\boldmath$B$}\in L^{*}$, let
\beqnarr
A+n<B+m & :\Lrarw & o(A)+n<o(B)+m
\label{eq:less}
\\
A\leq B & :\Lrarw & o(A)\leq o(B)
\nonumber
\\
A\simeq B & :\Lrarw & o(A)= o(B)
\nonumber
\\
\mbox{\boldmath$A$}<\mbox{\boldmath$B$} & :\Lrarw & 
\exi B\in \mbox{\boldmath$B$}\fal A\in\mbox{\boldmath$A$}(A<B)
\nonumber
\\
\mbox{\boldmath$A$}\leq\mbox{\boldmath$B$} & :\Lrarw & 
\fal A\in\mbox{\boldmath$A$}\exi B\in \mbox{\boldmath$B$}(A\leq B)
\nonumber
\eeqnarr
where $<,\leq$ in the RHS denote the relations in $O(\mu)$.
}
\edf



  
  
We are going to define \textit{moves} of hydras.
For a pair $(H,lb_{0})$
of a hydra $H$ and a finite set $lb_{0}$ of labels in $L$,
there are some possible moves $(H,lb_{0})\to_{\ell} (K,lb_{1})$ 
 depending on a number $\ell<\ome$.
The finite sets of labels may grow
in two cases (Production) in Definition \ref{df:move}.\ref{move:bracemu} and \ref{df:move}.\ref{move:braceB}.

 \bdf\label{df:move}{\rm (Moves)}\\
{\rm Let $(H,lb)$ be a pair of a hydra $H$ and a finite set $lb$ of labels in $L$, and $\ell<\ome$.
We define possible moves $(H,lb)\to_{\ell}(K,lb^{\prime})$.
}

 \benu
 \item\label{move:1}
{\rm (Necrosis)
$(H,lb)\to_{\ell}(0,lb)$ for $H\neq 0$,
$(H,lb)\to_{\ell}(1,lb)$ for $H\not\in\{0,1\}$
and
$(\{\mu\}(H),lb)\to_{\ell} (H,lb)$
.}
  

  \item\label{move:cdot*}
{\rm $(n\cdot *_{\mu},lb)\to_{\ell} ((n\cdot A)+n,lb)$ for $A\in lb$.}

{\rm $(n\cdot B,lb)\to_{\ell} \left((n\cdot A)+n,lb\right)$ for $lb\ni A< B$,
where 
$n=1+\cdots+1$.}

{\rm $(n\cdot *_{\ome},lb)\to_{\ell} (n\cdot m,lb)$ for $0<m\leq \ell$.}
      
  {\rm  $((n+1)\cdot (m+1),lb)\to_{\ell} (((n+1)\cdot m)+n,lb)$ for $n\geq 0$ and $m>0$. 
  $((n+1)\cdot 1,lb)\to_{\ell} (n,lb)$.}
   
\item\label{move:d}
 {\rm 
 $(d(H+1),lb)\to_{\ell}(d(H)\cdot k,lb)$, where $k\leq \ell+1$ and
 $d(H)\cdot k:=d(H)+\cdots+d(H)$ with $k$'s $d(H)$ for 
 $d\in\{\ome\}\cup\{D_{\mbox{\scriptsize \boldmath$C$}}, \vphi_{C+n}: \mbox{\boldmath$C$}\in L^{*}, C\in L, n\leq\ell\}$.}
     
\item\label{move:D2}
{\rm
$(D_{\mbox{\scriptsize \boldmath$C$}}(H+1),lb)\to_{\ell}(\vphi_{A+n}(D_{\mbox{\scriptsize \boldmath$C$}}(H)\cdot 2),lb)$ where $A\in lb$ with 
$A\leq \mbox{\boldmath$C$}$, and $n\leq\ell$.
}


\item\label{move:vphi2}
{\rm
$(\vphi_{C+n}(H+1),lb)\to_{\ell}(\vphi_{A+m}(\vphi_{C+n}(H)+\vphi_{\mbox{\scriptsize \boldmath$B$}}(H)),lb)$ and 
\\
$(\vphi_{C+n}(H+1),lb)\to_{\ell}(\vphi_{A+m}(\vphi_{\mbox{\scriptsize \boldmath$B$}}(H)+\vphi_{C+n}(H)),lb)$,
\\
where $A\in lb$, $m\leq \ell$, $A+m<C+n$, $\mbox{\boldmath$B$}\subset lb$ and $\mbox{\boldmath$B$}<C$.

}

 \item\label{move:brace1}
 {\rm $(\{*_{\mu}\}(H),lb)\to_{\ell} (\{A\}(H),lb)$.
 }
 

\item\label{df:move.brace2}
{\rm $(d(K+\{B\}(H)),lb)\to_{\ell} (\{B\}(d(K+H)\cdot 2),lb)$
where $R(B)$, and $d=\ome$ if $B=\mu$.
Otherwise 
$d\in\{\ome\}\cup\{\vphi_{A+n}: B\leq A, 0<n\leq\ell\}$.

For $B\in L$ with $B<\mu$ and $\mbox{\boldmath$C$}\in L^{*}$,
$(D_{\mbox{\scriptsize \boldmath$C$}}(K+\{B\}(H)),lb)\to_{\ell} (\{B\}(D_{\mbox{\scriptsize \boldmath$C$}}(K+H)\cdot 2),lb)$.
}

\item\label{move:bracemu}
{\rm (Production) For $A=d_{\mu}(K+\{B\}(H))$ with $lb\cup\{0\}\ni B<D(\mbox{\boldmath$C$};K+\{\mu\}(H))$, $n\leq\ell$,
\\
$(D(\mbox{\boldmath$C$};K+\{\mu\}(H)),lb)\to_{\ell} (\vphi_{A+n}(D(\mbox{\boldmath$C$}\cup\{A\};K+H)\cdot 2),lb\cup\{A\})$.



}

\item\label{move:braceB}
{\rm (Production) Let $R(B)$, $H\in\calh_{1}$ and 
$lb\cup\{0\}\ni C<B$.
Also $e(*)$ is a hydra with a hole $*$ generated from the hole $*$ by
applying $H(*)\mapsto K+H(*), \vphi_{C_{0}+m}(H(*)), \vphi_{\mbox{\scriptsize \boldmath$C$}}(H(*))$
for $lb\supset \{C_{0}\}\cup \mbox{\boldmath$C$}<B$.
Then for $A=d_{B}(\{C\}(e(H)))$ and $n\leq\ell$,
$
(e(\{B\}(H)),lb)\to_{\ell}
\left(\vphi_{A+n}(\vphi_{A}(e(H))\cdot 2), lb\cup\{A\}\right)
$.
}
\item\label{move:D}
{\rm If 
$(H,lb)\to_{\ell} (K,lb)$ for $H\in\calh_{0}$
and $\fal A\in lb(A<D_{\mbox{\scriptsize \boldmath$C$}}(H))$, then
$(D_{\mbox{\scriptsize \boldmath$C$}}(H),lb)\to_{\ell} (D_{\mbox{\scriptsize \boldmath$C$}}(K),lb)$.
}

\item\label{move:+}
{\rm If $(H,lb_{0})\to_{\ell} (K,lb_{1})$, then
$(d(H),lb_{0})\to_{\ell} (d(K),lb_{1})$ for
$d\in\{\ome\}\cup\{\vphi_{A+n},\vphi_{\mbox{\scriptsize \boldmath$A$}}: A\in L, n\leq\ell, \mbox{\boldmath$A$}\in L^{*}\}\cup\{H_{0}+: H_{0}\in\calh\}$.
}


 \eenu

{\rm
$(H,lb)\to_{\ell}^{*}(K,lb^{\prime})$ denotes the reflexive and transitive closure of the relation
$\to_{\ell}$.
}

\edf

Definition \ref{df:move}.\ref{move:d} means that
$(\ome(H+1),lb)\to_{\ell}(\ome(H)\cdot 2,lb)$, 
$(D_{\mbox{\scriptsize \boldmath$C$}}(H+1),lb)\to_{\ell} (D_{\mbox{\scriptsize \boldmath$C$}}(H)\cdot 2,lb)$ and
$(\vphi_{C+n}(H+1),lb)\to_{\ell}(\vphi_{C+n}(H)\cdot 2,lb)$.
Definition \ref{df:move}.\ref{move:+} means that
$(\ome(H),lb_{0})\to_{\ell} (\ome(K),lb_{1})$,
$(\vphi_{A+n}(H),lb_{0})\to_{\ell} (\vphi_{A+n}(K),lb_{1})$, 
$(\vphi_{\mbox{\scriptsize \boldmath$A$}}(H),lb_{0})\to_{\ell} 
(\vphi_{\mbox{\scriptsize \boldmath$A$}}(K),lb_{1})$ and
$(H_{0}+H,lb_{0})\to_{\ell} (H_{0}+K,lb_{1})$
if $(H,lb_{0})\to_{\ell} (K,lb_{1})$.

It is clear that both of the relations $(H,lb)\to_{\ell}(K,lb^{\prime})$ and $A<B$ elementary recursive
on hydras $H,K$, finite sets $lb,lb^{\prime}$, labels $A,B$ and numbers $\ell$.
Moreover when $(H,lb)\to_{\ell}(K,lb^{\prime})$, either $lb^{\prime}=lb$ or $lb^{\prime}=lb\cup\{A\}$
for a label $A$.
\\

Given a hydra $H_{0}$ and a finite set $lb_{0}$ of labels in $L$, 
a finitely branching tree $Tr(H_{0},lb_{0})=\{t\in {}^{<\ome}\ome: (H_{0}[t],lb_{0}[t]) \mbox{ {\rm is defied}}\}$
 is obtained as follows.
$H_{0}[\eps]:=H_{0}$ and $lb[\eps]:=lb_{0}$ for the empty sequence $\eps$.
$\{(H_{0}[t*(i)], lb[t*(i)])\}_{i}$ is the set of the pairs $(K,lb^{\prime})$ such that
$(H_{0}[t],lb[t])\to_{\ell}(K,lb^{\prime})$ for the length $\ell=|t|$ of $t\in{}^{<\ome}\ome$.
We see that the tree $Tr(H_{0},lb_{0})$ is elementary recursive.

\section{Provability}\label{sect:prov}

We show
 the following holds as long as $H_{0}[t*(k)]$ is defined for  $t\in{}^{<\ome}\ome$ and $lb_{0}[\eps]=lb_{0}$:
\beqn\label{eq:provable}
d_{\Ome}(o(H_{0}[t*(k)]\#o(lb_{0}[t*(k)]))<d_{\Ome}(o(H_{0}[t])\# o(lb_{0}[t]))
\eeqn
Then Theorem \ref{th:main}.\ref{th:main2} follows from \cite{odMahlo}.

\blem\label{lem:provable}
Let $(H,lb)\to_{\ell}(K,lb^{\prime})$, and $A$ be the label defined as follows.
If $lb^{\prime}=lb$, then let $A:=0$.
Otherwise $lb^{\prime}=lb\cup\{A\}$.

Then
$o(K)\#o(A)<o(H)$,
$\fal \sig[K_{\sig}(o(K))\leq K_{\sig}(o(H))\cup K_{\sig}(o(lb^{\prime}))]$
and $\fal\sig\fal \alp\in K_{\sig}(o(A))[\alp\in K_{\sig}(o(H))\cup K_{\sig}(o(lb)) \lor
\alp<d_{\sig}(o(H))]$, where
$K_{\sig}(o(lb^{\prime}))=\bigcup\{K_{\sig}(o(B)): B\in lb^{\prime}\}$.
\elem
\bprf
We show the lemma by main induction on the sum of the sizes
$\#(H)+\#(K)$ with subsidiary induction on the cardinality of the finite sets $lb$.

Consider the case in Definition \ref{df:move}.\ref{df:move.brace2}.
\\
First let $(d(K+\{B\}(H)),lb)\to_{\ell} (\{B\}(d(K+H)\cdot 2),lb)$
where $R(B)$, and $d=\ome$ if $B=\mu$.
Otherwise 
$d\in\{\ome\}\cup\{\vphi_{A+n}: B\leq A, 0<n\leq\ell\}$.
Let $\alp=o(K)$, $\bet=o(H)$ and $\sig=o(B)$.
First consider the case $d=\ome$.
Then $H,K\in\calh_{0}$ by Definition \ref{df:hydra}.\ref{df:hydra.3}, and 
$\eta:=o(d(K+\{B\}(H)))=\ome^{\alp\#\sig\#\bet}$, while
$\xi:=o(\{B\}(d(K+H)\cdot 2))=\sig\#\ome^{\alp\#\bet}\cdot 2$.
It is clear that $\xi<\eta$ and $\fal\tau(K_{\tau}\xi\leq K_{\tau}\eta)$.
Next let $d=D_{C}$ and $\sig<\mu$.
Then $H,K\in\calh_{0}$ by Definition \ref{df:hydra}.\ref{df:hydra.4}, and
$\eta:=o(d(K+\{B\}(H)))=d_{\mu}(\gam\#\alp\#\sig\#\bet\#1)$ with $\gam:=o(C)$, while
$\xi:=o(\{B\}(d(K+H)\cdot 2))=\vphi_{\sig}(d_{\mu}(\gam\#\alp\#\bet\#1)\cdot 2)$.
We see $\xi<\eta$ from $\sig<\eta$, and 
$\fal\tau<\mu(K_{\tau}\xi\subset K_{\tau}\{\gam,\alp,\sig,\bet\}\leq K_{\tau}\eta)$.
Also $K_{\mu}\xi=\{\sig,d_{\mu}(\gam\#\alp\#\bet\#1)\}<\{d_{\mu}(\gam\#\alp\#\sig\#\bet\#1)\}=K_{\mu}\eta$.
Finally let $d=\vphi_{A+n}$ with $B\leq A+n$ and 
$\rho=o(A)$.
Then $\eta:=o(d(K+\{B\}(H)))=\vphi_{\rho+n+1}(\alp\#\vphi_{\sig}(\bet))$ and 
$\xi:=o(\{B\}(d(K+H)\cdot 2))=\vphi_{\sig}(\vphi_{\rho+n+1}(\alp\#\bet)\cdot 2)$.
We have $\sig<\rho+n+1$.
We see $\xi<\eta$ from $\sig<\rho+n+1$ and $\bet<\vphi_{\sig}(\bet)$.
It is clear that $\fal\tau(K_{\tau}\xi\leq K_{\tau}\eta)$.

Second let for $B\in L$ with $B<\mu$ and $\mbox{\boldmath$C$}\in L^{*}$,
$(D_{\mbox{\scriptsize \boldmath$C$}}(K+\{B\}(H)),lb)\to_{\ell} 
(\{B\}(D_{\mbox{\scriptsize \boldmath$C$}\cup\{B\}}(K+H)\cdot 2),lb)$.
Let $\alp=o(K)$, $\bet=o(H)$, $\gam=o(\mbox{\boldmath$C$})$ and $\sig=o(B)$.
Then $o(D_{\mbox{\scriptsize \boldmath$C$}}(K+\{B\}(H)))=d_{\mu}(\gam\#\sig\#\alp\#\bet\#1)$, while
$o((\{B\}(D_{\mbox{\scriptsize \boldmath$C$}\cup\{B\}}(K+H)\cdot 2))=\vphi_{\sig}(d_{\mu}(\gam\#\sig\#\alp\#\bet))$.
It is clear that $\vphi_{\sig}(d_{\mu}(\gam\#\sig\#\alp\#\bet))<d_{\mu}(\gam\#\sig\#\alp\#\bet\#1)$
and $K_{\tau}\vphi_{\sig}(d_{\mu}(\gam\#\sig\#\alp\#\bet))=K_{\tau}(\gam,\sig,\alp,\bet)$ for $\tau<\mu$.

Next consider the case in Definition \ref{df:move}.\ref{move:bracemu}.
\\
$(D(\mbox{\boldmath$C$};K+\{\mu\}(H)),lb)\to_{\ell} (\vphi_{A+n}(D(\mbox{\boldmath$C$}\cup\{A\};K+H)\cdot 2),lb\cup\{A\})$, where
$n\leq\ell$, 
 $A=d_{\mu}(K+\{B\}(H))$with 
$lb\ni B<D(\mbox{\boldmath$C$};K+\{\mu\}(H))$.
Let $\alp=o(K)$, $\bet=o(H)$ for $K,H\in\calh_{0}$,
$\gam_{B}=o(B)$, $\gam=o(\mbox{\boldmath$C$})$.
Then $\del:=o(A)=d_{\mu}(\gam_{B}\#\alp\#\bet\#1)$,
$\xi:=o(\vphi_{A+n}(D(\mbox{\boldmath$C$}\cup\{A\};K+H)\cdot 2))=\vphi_{\del+n+1}(d_{\mu}(\gam\#\del\#\alp\#\bet)\cdot 2)$ and
$\eta:=o(D(\mbox{\boldmath$C$};K+\{\mu\}(H)))=d_{\mu}(\gam\#\mu\#\alp\#\bet\#1)$.
We have $\gam_{B}<\eta$.
Hence $\del<\eta$ and $\xi<\eta$.
On the other hand we have for $\tau<\mu$,
$K_{\tau}\xi\subset K_{\tau}\{\gam,\del,\alp,\bet\}=K_{\tau}\gam_{B}\cup K_{\tau}\eta$ with
$K_{\tau}\gam_{B}\subset K_{\tau}o(lb)$.
Next for $\tau=\mu$,
$K_{\mu}\xi=\{\del,d_{\mu}(\gam\#\del\#\alp\#\bet)\}<\{\eta\}=K_{\mu}\eta$.

Third consider the case in Definition \ref{df:move}.\ref{move:braceB}.
\\
$(e(\{B\}(H)),lb)\to_{\ell}\left(\vphi_{A+n}(\vphi_{A}(e(H))\cdot 2), lb\cup\{A\}\right)$,
where $R(B)$, $H\in\calh_{1}$, $n\leq\ell$ and $A=d_{B}(\{C\}(H))$ with
 $lb\cup\{0\}\ni C<B$.
$e(*)$ is a hydra with a hole $*$ generated from the hole $*$ by
applying $H(*)\mapsto K+H(*), \vphi_{C_{0}+m}(H(*)), \vphi_{\mbox{\scriptsize \boldmath$C$}}(H(*))$
for $lb\supset \{C_{0}\}\cup \mbox{\boldmath$C$}<B$.
 
Let $\sig=o(B)$, $\bet=o(H)$, and $\gam=o(C)$.
Also $\alp(*)=o(e(*))$ built from $+$ and $\vphi_{\kap}$ with $\kap<\sig$.
Then $\del:=o(A)=d_{\sig}(\sig\#\vphi_{\gam}(\alp(\bet)))$ and
$\xi:=o(\vphi_{A+n}(\vphi_{A}(e(H))\cdot 2))=
\vphi_{\del+n+1}(\vphi_{\del+1}(\alp(\bet))\cdot 2)$ and
$\eta:=o(e(\{B\}(H)))=\alp(\vphi_{\sig}(\bet))$.
From $\del<\sig$ and $\max\{\sig,\bet\}<\vphi_{\sig}(\bet)$ we see
$\xi\#\del<\eta$.
We have $K_{\tau}\{\xi,\del\}\subset
K_{\tau}\{\del,\bet,\alp(0)\}\subset K_{\tau}\del\cup K_{\tau}\eta$ for any $\tau$.
First let $\tau<\sig$. 
Then $K_{\tau}\del=K_{\tau}\{\sig,\gam,\bet,\alp(0)\}\subset K_{\tau}\{\eta,\gam\}$,
and $K_{\tau}\gam=K_{\tau}o(C)\subset K_{\tau}o(lb)$.
Next let $\tau=\sig$.
We have $\gam<\sig$.
Proposition 5.1.7 in \cite{ptMahlo} yields $K_{\tau}\gam\leq K_{\tau}\sig<d_{\tau}(\eta)$ for any $\tau>\sig$.
Hence we obtain $K_{\sig}\del=\{\del\}<d_{\sig}(\eta)$.
Third let $\sig<\tau$.
Proposition 5.1.7 in \cite{ptMahlo} yields
$K_{\tau}\del\leq K_{\tau}\sig<d_{\tau}(\eta)$.

Fourth consider the case in Definition \ref{df:move}.\ref{move:D}.
$(D_{\mbox{\scriptsize \boldmath$C$}}(H),lb)\to_{\ell} (D_{\mbox{\scriptsize \boldmath$C$}}(K),lb)$ follows from
$(H,lb)\to_{\ell} (K,lb)$, where
$\fal A\in lb(A<D_{\mbox{\scriptsize \boldmath$C$}}(H))$.
Let $\alp=o(H)$, $\bet=o(K)$ and $\gam=o(\mbox{\boldmath$C$})$.
By IH we have $\bet<\alp$ and $\fal\sig(K_{\sig}\bet\subset K_{\sig}(\{\alp\}\cup o(lb)$.
For $\sig<\mu$, this yields 
$K_{\sig}d_{\mu}(\gam\#\bet)=K_{\sig}\{\gam,\bet\}\subset K_{\sig}(\{\gam,\alp\}\cup o(lb))=K_{\sig}d_{\mu}(\gam\#\alp\#1)\cup K_{\sig}o(lb)$.
On the other hand we have 
$K_{\mu}o(lb)\leq o(lb)<d_{\mu}(\gam\#\alp\#1)=o(D(\mbox{\boldmath$C$};H))$.
Hence $K_{\mu}o(D(\mbox{\boldmath$C$};K))=\{o(D(\mbox{\boldmath$C$};K))\}<\{o(D(\mbox{\boldmath$C$};H))\}=K_{\mu}o(D(\mbox{\boldmath$C$};H))$.

Finally consider the case in Definition \ref{df:move}.\ref{move:+}.
$(d(H),lb_{0})\to_{\ell} (d(K),lb_{1})$ follows from
$(H,lb_{0})\to_{\ell} (K,lb_{1})$ for
$d\in\{\ome\}\cup\{\vphi_{A+n},\vphi_{\mbox{\scriptsize \boldmath$A$}}: A\in L, n<\ome, \mbox{\boldmath$A$}\in L^{*}\}\cup\{H_{0}+: H_{0}\in\calh\}$.
Let $\alp=o(H)$ and $\bet=o(K)$.
By IH we have $\bet\#\del<\alp$ for $\del=o(B)$ with $lb_{1}\subset lb_{0}\cup\{B\}$.
Thus $o(d(K))\#\del<o(d(H))$.
On the other hand we have for any $\sig$,
$K_{\sig}o(d(K))\subset K_{\sig}\{\bet,o(d(H))\}\leq K_{\sig}(\{o(d(H))\}\cup o(lb_{1}))$
by IH,
and $d_{\sig}\alp<d_{\sig}o(d(H))$.
\eprf
\\

The following Corollary \ref{cor:provable} shows (\ref{eq:provable}).

\bcor\label{cor:provable}
If $(H,lb)\to_{\ell}(K,lb^{\prime})$, then $d_{\Ome}\left(o(K)\#o(lb^{\prime})\right)< d_{\Ome}\left(o(H)\#o(lb)\right)$.
\ecor

\section{Unprovability}\label{sect:unprov}

In this section we introduce first
a theory $[\Pi^{0}_{1},\Pi^{0}_{1}]$-Fix for non-monotonic inductive definitions
of $[\Pi^{0}_{1},\Pi^{0}_{1}]$-operators in \cite{Richter-Aczel74}.
In \cite{wienpi3} it is shown that the 1-consistency of the set theory {\sf KPM} is reduced to one of
the theory $[\Pi^{0}_{1},\Pi^{0}_{1}]$-Fix.
For Theorem \ref{th:main}.\ref{th:main3} it suffices to show, over {\sf EA},
 the 1-consistency of  $[\Pi^{0}_{1},\Pi^{0}_{1}]$-Fix
assuming the fact that the hydra game eventually terminates.

\subsection{A theory $[\Pi^{0}_{1},\Pi^{0}_{1}]$-Fix}

In \cite{wienpi3} we show that the wellfoundedness 
is provable up to each ordinal diagram $\alp<\Ome$ in a theory $[\Pi^{0}_{1},\Pi^{0}_{1}]$-Fix for 
$[\Pi^{0}_{1},\Pi^{0}_{1}]$-non-monotonic inductive definitions in \cite{Richter-Aczel74}.

For a class of formulas $\Phi$,
the theory $\Phi$-Fix for non-monotonic inductive definitions are two-sorted: one sort $x$ for natural numbers and the other $a$ for \textit{ordinals}.
The binary predicate $x\in I^{a}$, then, denotes the $a$-th stage of inductive definition
by a fixed operator $\Gam:\calP(\ome)\to\calP(\ome)$, which is defined by a first order formula 
$\Gam(X,x)\in\Phi$ in the language of the first order arithmetic $\calL({\sf PA})$ with an extra unary predicate $X$.
The axioms of the theories are:
\benu
\item
 Axioms of {\sf PA} and equality axioms for either sort.
\item 
The defining axiom of $x\in I^{a}$: $x\in I^{a}\lrarw \exi b<a[x\in\Gam(I^{b})]$.
\item 
Closure axiom: $\Gam(I^{\infty})\subset I^{\infty}$ for $I^{\infty}:=\{x:\exi a(x\in I^{a})\}$.
\item
Axioms for the well ordering $<$ on ordinals:
 \benu
 \item
 $<$ is a \textit{linear} ordering:
  \benu
  \item $<$ is irreflexive and transitive.
  \item (trichotomy) 
  \beqn\label{eq:trichotomy}
  x<y\lor x=y \lor y<x
  \eeqn
  \eenu
 \item
 transfinite induction schema for any formula $F$:
 \\
$\fal a[\fal b<aF(b)\rarw F(a)]\rarw \fal a F(a)$.
 \eenu
\eenu

\subsection{Hydras associated with proofs}

In what follows assume that $[\Pi^{0}_{1},\Pi^{0}_{1}]$-Fix is 1-inconsistent. 
This means that there exists a true $\Pi^{0}_{1}$-sentence 
$\fal x B(x)$ with a quantifier-free $B$ such that $[\Pi^{0}_{1},\Pi^{0}_{1}]$-Fix+$\fal x B(x)$  is inconsistent.

Let $P_{0}$ be a proof in $[\Pi^{0}_{1},\Pi^{0}_{1}]$-Fix+$\fal x B(x)$ of a contradiction. 
(Proofs are specified later.)
We associate a hydra $H_{0}=\Ome(P_{0})$ to $P_{0}$, and 
define a rewriting  step $r: P\mapsto r(P)$ on proofs $P$ in $[\Pi^{0}_{1},\Pi^{0}_{1}]$-Fix+$\fal x B(x)$.
For  each $P_{\ell}=r^{(\ell)}(P_{0})$, associate a hydra $H[\ell]=\Ome(P_{\ell})$ again so that
$\{H[\ell]\}_{\ell}$ is a path through the tree $Tr(H_{0},\emptyset)$.
$P_{0}$ tells the hydras which way to proceed.
Namely $H[\ell+1]$ is one of possible moves for  the hydra $H[\ell]$, i.e.,
$(H[\ell],lb[\ell])\to_{\ell}(H[\ell+1],lb[\ell+1])$.
Assuming $P_{0}$ is a proof in $[\Pi^{0}_{1},\Pi^{0}_{1}]$-Fix+$\fal x B(x)$ of a contradiction,
we see that the path is infinite, i.e., 
the hydra game $\{H[\ell]\}_{\ell}$ goes forever.
Moreover all of these are done in {\sf EA}.
\\








$\calL^{2}$ denotes the class of lower elementary recursive functions in \cite{Sk}.
The class of functions containing the zero, successor, projection and modified subtraction functions 
and which is closed under composition and summation of functions. 
$\calL^{2}_{*}$ denotes the class of lower elementary recursive relations. 
The arithmetical part of the language $\calL$ of $[\Pi^{0}_{1},\Pi^{0}_{1}]$-Fix
is chosen to consist of predicate constants for lower elementary recursive relations $R\in \calL^{2}_{*}$. 

The \textit{language} $\calL$ of $[\Pi^{0}_{1},\Pi^{0}_{1}]$-Fix consists of 
\benu
\item
two sorts of variables, one for (natural) numbers $\Natural$ and the other for ordinals, i.e., 
stages $\mathcal{O}$ of inductive definitions.
$x,y,\ldots$ are variables for natural numbers, and $a,b,\ldots$ are variables ranging over
the domain of a well ordering $<$,

\item
two binary predicate symbols $a=^{\mathcal{O}}b$ and $a< b$, and their negations $a\neq b$ and $a\not< b$
on $\mathcal{O}$,

\item 
function constants $0^{\Natural}$ and $x^{\prime}$ (successor) on $\Natural$, 

\item 
arithmetic predicate constants on $\Natural$ for lower elementary recursive relations $R\in \calL^{2}_{*}$ and their negations $\lnot R$,

\item 
the binary predicate symbol $I(a,x)$ and its negation $\lnot I(a,x)$ 
denoting the stages
$I^{a}=\{x\in\ome: I(a,x)\}$ of a fixed $[\Pi^{0}_{1},\Pi^{0}_{1}]$-formula 
\[
\cala(X,x)\equiv
\fal y\, \calb_{0}(X,x,y)\lor
 [\fal x\{\fal y\, \calb_{0}(X,x,y)\to x\in X\}\land \fal z\, \calb_{1}(X,x,z)]
\]
where $\calb_{i}$ is a bounded formula in the arithmetic language $\calL^{2}_{*}\cup\{0^{\Natural},\prime\}$ 
with a unary predicate $X$ for $i=0,1$, and

\item 
logical connectives $\land, \lor, \fal, \exi$.
\eenu

The \textit{negation} $\lnot \vphi$ of a formula $\vphi$ is defined by using de Morgan's law and the elimination of double negations. 
A prime formula $R(t_1,\ldots,t_n)$ or its negation $\lnot R(t_1,\ldots,t_n)$ 
with an arithmetic predicate $R$ is an \textit{a.p.f.}(arithmetic prime formula),
and a prime formula $t=s, t< s$ for stage terms $t,s$ and their negations are
\textit{s.p.f.}(stage prime formula).

There are four kinds of quantifications, \textit{bounded} number quantifiers $\exi x\leq t, \fal x\leq t$,
\textit{unbounded} number quantifiers $\exi x, \fal x$,
\textit{bounded} stage quantifiers $\exi a< b, \fal a< b$ and \textit{unbounded} stage quantifiers $\exi a,\fal a$.
A formula is said to be \textit{unbounded} if it contains an unbounded \textit{stage} quantifier.

The \textit{axioms} in $[\Pi^{0}_{1},\Pi^{0}_{1}]$-Fix are axioms for function and arithmetic predicate constants, 
the axioms 
for the linear ordering $<$,
the induction axioms $(VJ)$, $(TJ)$, the defining axiom $(I)$ of stages,
and the closure axiom $(Cl)$: for arbitrary formula $F$,\bdes
\item[$(VJ)$] 
$F(0) \land \fal x(F(x)\to F(x^{\prime}))\, \to \, \fal xF(x)$.

\item[$(TJ)$]
$\fal a[\fal b< a\, F(b)\to F(a)]\, \to \, \fal a F(a)$.

\item[$(I)$] 
$\fal x\fal a[x\in I^{a}\lrarw (\cala(I^{< a},x)\lor x\in I^{< a})]$,
where $(x\in I^{a}):\Lrarw I(a,x)$ and $(x\in I^{< a}):\Lrarw \exi b< a(x\in I^{b})$
with $I^{b}=\{x: I(b,x)\}$.

\item[$(Cl)$]
$\cala(I^{<\infty})\subset I^{<\infty}$, i.e.,
$\fal x(\cala(I^{<\infty},x) \to x\in I^{<\infty})$, where
$(x\in I^{<\infty}):\Lrarw \exi a(x\in I^{a})$.
This is equivalent to
$\fal y\, \calb_{0}(I^{<\infty},y)\subset I^{<\infty}
\land
  \fal z\, \calb_{1}(I^{<\infty},z) \subset I^{<\infty}
 $, where
$\fal y\, \calb_{0}(I^{<\infty},y)=\{x: \fal y\, \calb_{0}(I^{<\infty},x,y)\}$
and
$ \fal z\, \calb_{1}(I^{<\infty},z)=\{x: \fal z\, \calb_{1}(I^{<\infty},x,z)\}$.
\edes

Let us extend the language $\calL$ to $\calL_{H}$ by adding a unary predicate $R(a)$ of stage sort,
and individual constants $A$ denoting 
labels $A\in L$ and a constant $0^{\mathcal{O}}$ for the hydra $0\in\calh$. 
By definition these constants $A$ is of stage sort.
$A=^{\mathcal{O}}B$ is defined to be \textit{true} iff $o(A)=o(B)$,
and $A< B$ is true if $o(A)<o(B)$.
$lb(\vphi)$ denotes the set of stage constants occurring in the formula $\vphi$.
On the other side $R(t)$ and $\lnot R(t)$ ate s.p.f's, and $R(A)$ is defined to be \textit{true} iff $A=d_{\mu}(h)$
for some $h$.
$R(A)$ is intended to denote the fact that $A$ is recursively regular.
Then the axiom $(Cl)$, $\cala(I^{<\infty})\subset I^{<\infty}$ is proved from the following axioms.
\bdes
\item[$(Cl.0)$]
$\fal x\{\fal y\, \calb_{0}(I^{<\infty},x,y) \to\exi a[R(a)\land \fal y\,\calb_{0}(I^{< a},x,y)]\}$ and
\\
$\fal a\left(R(a)\to \fal x\{\fal y\, \calb_{0}(I^{< a},x,y)\to \exi b< a[\fal y\,\calb_{0}(I^{< b},x,y)\}\right)$

\item[$(Cl.1)$]
$\fal z\{\fal y\, \calb_{1}(I^{<\infty},z,y)\to \exi a[R(a)\land  \fal y\, \calb_{1}(I^{< a},z,y)]\}$.
\edes
In what follows by a formula we mean a formula in $\calL_{H}$.

A formula is said to be an $\exi$\textit{-formula}  if it is either an a.p.f. or a s.p.f. or a formula in one of the following shapes; $\vphi\lor\psi$, $\exi x\leq t\vphi$, $\exi x\vphi$, $\exi a< b\,\vphi$, $\exi a \vphi$ or $t\in I^{a}$.
A formula is a $\fal$\textit{-formula} if its negation is an $\exi$-formula. 
If a formula is an $\exi$-formula and simultaneously a $\fal$-formula,
then it is either an a.p.f. or a s.p.f.

$[\Pi^{0}_{1},\Pi^{0}_{1}]$-Fix is formulated in one sided sequent calculus. 
Finite sets of formulae are called a \textit{sequents}. 
Sequents are denoted by $\Gam, \, \Del, \: etc.$

\bdf\label{df:sequentcalculus}
{\rm
{\bf Axioms} in $[\Pi^{0}_{1},\Pi^{0}_{1}]$-Fix are:
\bdes
\item[logical axioms] $\Gam,\lnot \vphi,\vphi$ \\
where $\vphi$ is an a.p.f. or a s.p.f. or a formula of the shape $t\in I^{s}$.


\item[arithmetical axioms] 
 \benu
 \item 
 $\Gam,\Del_{R}$ 
 \\
where $\Del_{R}$ consists of a.p.f'.s and corresponds to the definition of a lower elementary relation $R$.
 \item 
 $\Gam,\vphi$ for a true closed a.p.f. $\vphi$.
 \item 
 $\Gam,\Del_0$ 
 \\
where there exists a sequent $\Del_1$ so that $\Del=\Del_0\cup\Del_1$ is an instance of a defining axiom for $R$ in $1$ and $\Del_1$ consists solely of false closed a.p.f.'s.

Any true closed a.p.f. in an arithmetical axiom is said to be a {\it principal formula\/} of the axiom.
 \eenu
\item[stage prime axioms] 
 \benu
 \item
 $\Gam,t_{0}\not< t_{0}$, $\Gam, t_{0}\not< t_{1}, t_{1}\not< t_{2},t_{0}< t_{2}$, and 
 $\Gam,t_{0}< t_{1},t_{0}=t_{1},t_{1}< t_{0}$
 for terms $t_{0},t_{1},t_{2}$ of stage sort.
 \item
$\Gam,\vphi$ for a true closed s.p.f. $\vphi$.
\item
$\Gam, t\not\in I^{0^{\mathcal{O}}}$.

Each true closed s.p.f. in a stage prime axiom is said to be a {\it principal formula\/} of the axiom.

 \eenu
\edes

 Observe that the relation `a sequent $\Gam$ is an axiom in $[\Pi^{0}_{1},\Pi^{0}_{1}]$-Fix' is elementary recursive and hence so is the relation `$P$ is a proof in $[\Pi^{0}_{1},\Pi^{0}_{1}]$-Fix' with the inference rules defined below. 
\\

\noindent
{\bf Inference rules} in $[\Pi^{0}_{1},\Pi^{0}_{1}]$-Fix are: \\
$(\land),(\lor), (b\fal)^{\Natural},(b\exi)^{\Natural}, (\fal)^{\Natural},(\exi)^{\Natural}, (b\fal)^{\mathcal{O}},(b\exi)^{\mathcal{O}},(\fal)^{\mathcal{O}},(\exi)^{\mathcal{O}},(I),(\lnot I)$,
$(cut) $ and $(VJ)$, $(TJ)$, $(Cl)$.
\benu
\item 
\textit{Basic rules}
$(\land),(\lor), (b\fal)^{\Natural},(b\exi)^{\Natural}, (\fal)^{\Natural},(\exi)^{\Natural}, (b\fal)^{\mathcal{O}},(b\exi)^{\mathcal{O}},(\fal)^{\mathcal{O}},(\exi)^{\mathcal{O}}, (I),(\lnot I)$: 
In these rules the principal formula is contained in the upper sequent. 
For example
\[
\infer[(\lor)]{\vphi_{0}\lor\vphi_{1},\Gam}{\vphi_{0}\lor\vphi_{1},\vphi_{i},\Gam}
\:
\infer[(b\exi)^{\Natural}]{\exi x\leq t\, \vphi(x),\Gam}{\exi x\leq t\, \vphi(x),  u\leq t\land \vphi(u),\Gam}
\:
\infer[(\exi)^{\Natural}]{\exi x\, \vphi(x),\Gam}{\exi x\, \vphi(x),  \vphi(u),\Gam}
\]
where $i=0,1$, $u$ is a number term.
The \textit{minor formula} of these rules are defined to be
the formula $\vphi_{i}$ in $(\lor)$, and
$\vphi(u)$ in $(b\exi)^{\Natural}, (\exi)^{\Natural}$, resp.
The term $u$ in $(b\exi)^{\Natural},(\exi)^{\Natural}$ is the \textit{witnessing term} of the rules.

\[
\infer[(b\exi)^{\mathcal{O}}]{\exi a< t \,\vphi(a),\Gam}{\exi a< t\, \vphi(a), s< t\land \vphi(s),\Gam}
\msfiv
\infer[(\exi)^{\mathcal{O}}]{\exi a\, \vphi(a),\Gam}{\exi a\, \vphi(a),  \vphi(s),\Gam}
\]
where $s,t$ are stage terms.
The {\it minor formula\/} of these rules are defined to be
the formula $\vphi(s)$ both in $(b\exi)^{\mathcal{O}}$ and in $(\exi)^{\mathcal{O}}$.

For a number term $t$ and a stage term $s$,
\[
\infer[(I)]{t\in I^{s},\Gam}
{
t\in I^{s},\cala(I^{< s},t),\Gam
}
\msfiv
\infer[(I)]{t\in I^{s},\Gam}
{
t\in I^{s},t\in I^{< s},\Gam
}
\]
$\cala(I^{< s},t)$ and $t\in I^{< s}$ are the \textit{minor formula} of the rules $(I)$.
\[
\infer[(\lnot I)]{\Gam,t\not\in I^{s}}
{
\Gam,t\not\in I^{s},\lnot\cala(I^{< s},t)
&
\Gam,t\not\in I^{s},t\not\in I^{< s}
}
\]

\item 
In the rule $(cut)$ 
\[
\infer[(cut)]{\Gam,\Del}{\Gam,\lnot \vphi & \vphi,\Del}
\]
the \textit{cut formula} $\vphi$ is an $\exi$-formula.


\item 
\[
\infer[(VJ)]{\Gam}{\Gam, \vphi(0) & \lnot \vphi(x),\Gam, \vphi(x') & \lnot \vphi(t),\Gam}
\]
for any $\fal$\textit{-formulae} $\vphi$ and number terms $t$, where $x$ is the \textit{eigenvariable}.

$t$ is said to be the \textit{induction term}
of the $(VJ)$.

\item
\[
\infer[(TJ)]{t\not< s, \Gam}
{
\lnot \vphi(t),\Gam
&
\Gam,\lnot\fal b< a\, \vphi(b), \vphi(a)
}
\]
for any $\fal$\textit{-formulae} $\vphi$ and any stage terms $t,s$, where $a$ is the 
\textit{eigenvariable}.

$\vphi(a)$ is said to be the \textit{induction formula} and $s$ the \textit{induction term}
of the $(TJ)$.

\item
For  an \textit{eigenvariable} $a$,
\[
\infer[(Cl. R)]{\lnot R(s), \Gam}
{
\Gam, \fal y\, \calb_{0}(I^{< s},t,y)
&
a\not< s,  \lnot \fal y\, \calb_{0}(I^{< a},t,y),\Gam
}
\]
where $s$ denotes either a stage variable or a constant $A\in L$.
\[
\infer[(Cl.\mu)]{\Gam}
{
\Gam, \fal y\, \calb_{i}(I^{<\mu},t,y)
&
\lnot R(a), \lnot \fal y\, \calb_{i}(I^{< a},t,y),\Gam
}
\]
where $i=0,1$.

For a stage variable or a stage constant $s$, let
$R_{s}(a):\equiv(0^{\mathcal{O}}=0^{\mathcal{O}})$.
$(R(\mu)):\equiv(0^{\mathcal{O}}=0^{\mathcal{O}})$,
 and $(R_{\mu}(a)):\equiv(R(a))$.
These two rules are then unified to the following rule:
\[
\infer[(Cl.s)]{\lnot R(s), \Gam}
{
\Gam, \fal y\, \calb_{i}(I^{< s},t,y)
&
\lnot R_{s}(a),a\not< s,  \lnot \fal y\, \calb_{i}(I^{< a},t,y),\Gam
}
\]
where $s$ denotes either a stage variable or a constant $A\in L\cup\{\mu\}$, and $i=0$ when $s\neq\mu$.
\eenu
}
\edf

\bdf\label{df:tautology}
{\rm For each formula} $\vphi$ {\rm and sequent} $\Gam$, 
{\rm let} $P_{\Gam,\vphi}$ {\rm denote a canonically constructed proof of} $\Gam,\lnot\vphi,\vphi$
{\rm using logical axioms and rules} 
$(\lor), (\land)$, $(b\exi)^{\Natural},(b\fal)^{\Natural}$, 
$(\exi)^{\Natural},(\fal)^{\Natural}$, $(b\exi)^{\mathcal{O}},(b\fal)^{\mathcal{O}}$, 
$(\exi)^{\mathcal{O}}, (\fal)^{\mathcal{O}}$.

\benu
\item
{\rm If} $\vphi$ {\rm is an a.p.f. or a s.p.f. or a formula of the shape} $t\in I^{s}$, 
{\rm then} $P_{\Gam,\vphi}$ {\rm denotes the logical axiom}
$\Gam,\lnot\vphi,\vphi$.

\item
{\rm If} $\vphi\equiv(\tht_{0}\lor\tht_{1})$, {\rm then for} $\Del=\Gam,\lnot\vphi,\vphi$
\[
P_{\Gam,\vphi}=\left.
\begin{array}{c}
\infer[(\land)]{\Del}
{
 \infer[(\lor)]{\Del,\lnot\tht_{0}}
 {
  \infer*[P_{\Del,\tht_{0}}]{\Del,\lnot\tht_{0},\tht_{0}}{}
  }
 &
  \infer[(\lor)]{\Del,\lnot\tht_{1}}
 {
  \infer*[P_{\Del,\tht_{1}}]{\Del,\lnot\tht_{1},\tht_{1}}{}
  }
}
\end{array}
\right.
\]

\item
{\rm If} $\vphi\equiv(\exi a< t\,\tht(a))$, 
{\rm then for} $\Del=\Gam\cup\{\lnot\vphi,\vphi,a\not< t\lor\lnot\tht(a),a< t\land\tht(a)\}$,

\[
P_{\Gam,\vphi}=\left.
\begin{array}{c}
\infer[(b\fal)^{\mathcal{O}}]{\Gam,\lnot\vphi,\vphi}
{
 \infer[(b\exi)^{\mathcal{O}}]{\Gam,\lnot\vphi,\vphi,a\not< t\lor \lnot\tht(a)}
 {
  \infer[(\land)]{\Gam,\lnot\vphi,\vphi,a\not< t\lor\lnot\tht(a),a< t\land\tht(a)}
  {
   \infer[(\lor)]{\Del,a< t}
   {
    \Del,a\not< t,a< t
    }
  &
   \infer[(\lor)]{\Del,\tht(a)}
   {
    \infer*[P_{\Del,\tht}]{\Del,\lnot\tht(a),\tht(a)}{}
    }
   }
  }
}
\end{array}
\right.
\]

\item
{\rm And similarly for the cases} $\vphi\equiv(\exi x<t\,\tht(x)), (\exi x\,\tht(x)), (\exi a\,\tht(a))$.

\eenu
\edf

\bdf\label{df:rank}
{\rm 
The \textit{rank} $\rk(\vphi)\in \{A+n: A\in lb(\vphi)\cup\{0,\mu\}, n<\ome\}$
and the \textit{label complexity} $lq(\vphi)\in lb(\vphi)\cup\{0,\mu\}$
of a formula $\vphi$ in $\calL_{H}$ are defined recursively.
Let $Q\in\{\fal ,\exi\}$.
}
\benu

\item
$\rk(\vphi)=lq(\vphi)=0$ {\rm for an a.p.f. or a s.p.f.} $\vphi$. 

\item 
{\rm
$\rk(Qx\leq t \, \vphi)=\rk(Qx \, \vphi)=\rk(\vphi)+1$  and
$lq(Qx\leq t \, \vphi)=lq(\vphi)$ for the number variable $x$.
}

\item 
{\rm For 
$\circ\in\{\land,\lor\}$,
$\rk(\vphi_0\circ \vphi_1)=\max\{\rk(\vphi_{0}),\rk(\vphi_{1})\}+1$
and $lq(\vphi_0\circ \vphi_1)=\max\{lq(\vphi_{0}),lq(\vphi_{1})\}$.
}

\item
{\rm 
$\rk(Q a\,\vphi(a))=\mu$ if $\vphi(a)$ is bounded. Otherwise
$\rk(Q a\,\vphi(a))=\rk(\vphi(0^{\mathcal{O}}))+1$.
$lq(Q a\,\vphi(a))=\mu$.
}

\item
{\rm 
For contants $\mu\neq B\in L$,
$\rk(Q a< B\,\vphi(a))=\max\{\rk(\vphi(0^{\mathcal{O}}))+1,B\}$.
$lq(Q a< B\,\vphi(a))=\max\{lq(\vphi),B\}$ for $B\in L\cup\{0\}$.
}

\item
$\rk(Q a< 0^{\mathcal{O}}\,\vphi(a))=\rk(\vphi(0^{\mathcal{O}}))$.

\item
{\rm 
For variables $b$,
$\rk(Q a< b\,\vphi)=\mu$
if $Q a< B\,\vphi(a)$ is bounded. Otherwise
$\rk(Q a< b\,\vphi)=\rk(\vphi(0^{\mathcal{O}}))+1$.
$lq(Q a< b\,\vphi)=lq(\vphi)$.
}

\item
{\rm For $\mu\neq B\in L$,
$\rk(t\in I^{B})=\rk(t\not\in I^{B})=B+(d_{\cala}+1)$,
{\rm where} $d_{\cala}=\rk(\cala(I^{0^{\mathcal{O}}},t))$ 
 denotes the depth of (number) quantifiers and propositional connectives
 $\land,\lor$ {\rm in} $\cala(X,t)$.
$lq(t\in I^{B})=lq(t\not\in I^{B})=B$ for $B\in L\cup\{0\}$.
}

\item
$\rk(t\in I^{0^{\mathcal{O}}})=\rk(t\not\in I^{0^{\mathcal{O}}})=0$.

\item
{\rm 
For variables $b$,
$\rk(t\in I^{b})=\rk(t\not\in I^{b})=\mu$
and
$lq(t\in I^{b})=lq(t\not\in I^{b})=0$.
}

\eenu
\edf

Observe that $\rk(t\in I^{<\infty})=\rk(\exi a(t\in I^{a}))=\mu$,
while
$\rk(t\in I^{< A})=\rk(\exi a< A(t\in I^{a}))=\max\{\rk(t\in I^{0^{\mathcal{O}}})+1,A\}=A$ for $A\in L$.

\blem\label{lem:rank}
For any constants $\mu\neq A,B\in L$, the following hold.
\benu
\item\label{lem:rank0}

Let $\vphi$ be a closed formula.
Then
$\rk(\vphi)<\mu$ iff $\vphi$ is bounded, and
$\rk(\vphi)=lq(\vphi)+n$ for an $n<\ome$.

\item\label{lem:rank10}
$\rk(\vphi)<\mu$ iff
$\vphi$ is bounded, and
there occurs no subfoumulas $Qa< b\,\tht$, $t\in I^{b}$, $t\not\in I^{b}$
with variables $b$ in $\vphi$.

\item\label{lem:rank1}
$\rk(\lnot\vphi)=\rk(\vphi)$.

\item\label{lem:rank15}
For each formula $\vphi$, there exists a label $A\in lb(\vphi)\cup\{0,\mu\}$ such that
$\rk(\vphi)\leq A+\max\{d_{\cala}+1,n\}$, where $n$ denotes the number of occurrences of
logical connectives $\land,\lor,\fal,\exi$ in $\vphi$.

\item\label{lem:rank2}
$\rk(\vphi_{i})<\rk(\vphi_{0}\lor\vphi_{1})$ for $i=0,1$.


\item\label{lem:rank4}
$\rk(\vphi(\bar{n}))<\rk(\exi x \vphi(x))$ for the $n$-th numeral $\bar{n}$.




\item\label{lem:rank7}
$\rk(\vphi(A))<\rk(\exi a\, \vphi(a))$
if $\vphi(A)$ is closed.


\item\label{lem:rank8}
Assume that $A<B$. Then
$\rk(\vphi(A))<\rk(\exi a< B\, \vphi(a))$.

\item\label{lem:rank9}
$\rk(\cala(I^{< A},\bar{n}))<\rk(\bar{n}\in I^{A})$.

\eenu
\elem
\bprf
Lemma \ref{lem:rank}.\ref{lem:rank7} and \ref{lem:rank}.\ref{lem:rank8} follow from 
the facts that
$\rk(\vphi(0^{\mathcal{O}}))=\rk(\vphi(A))$ for unbounded $\vphi$,
and for bounded $\vphi$,
$\rk(\vphi(A))\in\{\rk(\vphi(0^{\mathcal{O}}))\}\cup\{A+n:n<\ome\}$.

For Lemma \ref{lem:rank}.\ref{lem:rank9} first observe that
$\rk(t\in I^{< A})=A$.
This yields $\rk(\cala(I^{< A},\bar{n}))=A+ d_{\cala}<\rk(\bar{n}\in I^{A})$.
\eprf

\bdf
{\rm 
We write $Qa<\mu$ for unbounded stage quantifier $Qa$.
For stage constants $A$ and formulas $\vphi$,
$\vphi^{A}$  denotes the result of restricting any unbounded stage quantifiers
$Q a<\mu$ to $Q a< A$ in $\vphi$.

$\Gam^{A}:=\{\vphi^{A}:\vphi\in\Gam\}$ for sequents $\Gam$.

For example $(\fal y\, \calb_{i}(I^{<\mu},x,y))^{A}\equiv(\fal y\, \calb_{i}(\{z:\exi a(z\in I^{a})\},x,y))^{A}\equiv(\fal y\, \calb_{i}(\{z: \exi a< A(z\in I^{a})\},x,y)\equiv(\fal y\, \calb_{i}(I^{< A},x,y))$.


}
\edf

The following definition is needed to handle bounded number quantifiers and propositional connectives, 
cf.\, subsections \ref{subsec:inversion} and \ref{subsec:topbnd}.

\bdf\label{df:resolvent}
{\rm
\textit{Resolvents} of a (closed) formula $\vphi$ are defined recursively as follows.
\benu
\item
$\Del=\{\vphi\}$ is a resolvent of $\vphi$.
\item
There is a resolvent $\Del_{1}\cup\{\tht_{0}\lor\tht_{1}\}$ such that $\Del=\Del_{1}\cup\{\tht_{0},\tht_{1}\}$.
\item
There is a resolvent $\Del_{1}\cup\{\tht_{0}\land\tht_{1}\}$ such that $\Del=\Del_{1}\cup\{\tht_{i}\}$ for an $i=0,1$.
\item
There is a resolvent $\Del_{1}\cup\{\exi x\leq m\tht(x)\}$ such that $\Del=\Del_{1}\cup\{\tht(k):k\leq m\}$.
\item
There is a resolvent $\Del_{1}\cup\{\fal x\leq m\tht(x)\}$ such that $\Del=\Del_{1}\cup\{\tht(k)\}$ for a $k\leq m$.
\eenu
}
\edf

Let $A\in L\cup\{\mu\}$.
A bounded formula $\tht$ is a \textit{$\Del^{A}$-formula} if each stage constant $C$ occurring in $\tht$ 
is $C<A$.
A formula $\vphi$ is a \textit{$\Sig^{A}$-formula} if either $\vphi$ is $\Del^{A}$ or $\vphi\equiv(\exi a<A\tht(a))$
with a $\Del^{A}$-formula $\tht$.
A \textit{$\Pi^{A}$-formula} is defined to be the dual of a $\Sig^{A}$-formula.

\bdf The system $[\Pi^{0}_{1},\Pi^{0}_{1}]\mbox{{\rm -Fix}} +\fal xB(x)$\\
{\rm Let} $\fal x B(x)$ {\rm denote a fixed} true $\Pi^0_1${\rm -sentence with an a.p.f.}\,$B${\rm . The system} 
$[\Pi^{0}_{1},\Pi^{0}_{1}]\mbox{{\rm -Fix}} +\fal xB(x)$ {\rm is obtained from} $[\Pi^{0}_{1},\Pi^{0}_{1}]\mbox{{\rm -Fix}}$ {\rm by adding the axioms} 
\[
(B)\; \Gam,B(t)
\]
{\rm for arbitrary terms} $t$ {\rm of number sort, and five inference rules;}
 {\rm the \textit{padding} rule ${}_{H_{0}}(pad)_{H_{1}}$, 
  the \textit{resolvent} rule $(res)_{H}$, the \textit{rank} rule $(rank)_{\mbox{\scriptsize \boldmath$C$}}$
  the \textit{height} rule $(h)$,
and the \textit{collapsing} rule} $(c)^{A}_{A(H)}$
{\small
 \[
 \infer[{}_{H_{0}}(pad)_{H_{1}}]{\Gam,\Del}{\Gam}
\hspace{3mm}
 \infer[(res)_{H}]{\bigcup_{k\leq n}\Pi_{k},\Del_{0}}
{
 t\leq \bar{n}\land\vphi(t),\bigcup_{k\leq n}\Pi_{k},\Del_{0}
 }
\hspace{3mm}
  \infer[(rank)_{\mbox{\scriptsize \boldmath$C$}}, (h)]{\Gam}{\Gam}
\hspace{3mm}
\infer[(c)^{A}_{d_{A}(H)}]{\Gam^{d_{A}(H)}}{\Gam^{A}}
\]
}
{\rm where $R(A)$ and} $\Gam$ {\rm denotes 
a finite set of closed subformulas of} $\Pi^{0}_{1}${\rm -formulae}
$\fal y\,\calb_{i}(I^{< \infty},n,y)$ {\rm with numerals $n$.
Each formula in $\Gam^{A}$ is obtained from $\Sig^{A}$-formulas, $\Pi^{A}$-formulas by
propositional connectives $\lor,\land$ and bounded number quantifications $\exi x\leq t,\fal x\leq t$.
In $(res)$ each $\Pi_{k}$ is a resolvent of the formula $\vphi(k)$, and $t\equiv y, y^{\prime}$ for a variable $y$.
The formula $t\leq \bar{n}\land\vphi(t)$ is the \textit{minor formula} of the $(res)$.}
\edf

In $(res)$,
$\vphi$ is a subformula of one of $\calb_{i}(I^{< B},n,m)$ and $\lnot\calb_{i}(I^{<  B},n,m)$
for some numerals $n,m$.

A \textit{proof} in the system $[\Pi^{0}_{1},\Pi^{0}_{1}]\mbox{{\rm -Fix}} +\fal xB(x)$ is a finite labelled tree
of sequents which is locally correct with respect to the axioms and inference rules in Definition 
\ref{df:sequentcalculus}.

\bdf\label{df:0prf}
{\rm Let} $P$ {\rm be a proof.}
\benu
\item 
{\rm For finite sequences ${\tt t}, {\tt s}\in {}^{<\ome}\ome$ of natural numbers,
${\tt t}\subset_{e}{\tt s}$ iff {\tt t} is an initial segment of {\tt s} in the sense that
${\tt s}={\tt t}*{\tt u}$ for a ${\tt u}{}^{<\ome}\ome$.
}
\item
$\mbox{{\rm Tr}}(P)\subset{}^{<\ome}\ome$ {\rm denotes the} underlying tree {\rm of} $P$,
{\rm where the endsequent corresponds to the root} $\eps$
{\rm (the} empty sequence{\rm ), and if a lowersequent} $\Gam$ {\rm of a rule} ${\tt t}*(0):J$ 
{\rm corresponds to a node {\tt t}, then 
 its uppersequents} $\Lam_{0},\ldots,\Lam_{n}$ {\rm correspond to} ${\tt t}*(0,0),\ldots,{\tt t}*(0,n)${\rm , resp.}
\[
\infer[{\tt t}*(0): J]{{\tt t}:\Gam}{{\tt t}*(0,0):\Lam_0 & \cdots & {\tt t}*(0,n):\Lam_n}
\]

\item 
{\rm For a node} ${\tt t}$, ${\tt t}:\Gam$ {\rm designates that 
the sequent $\Gam$ is situated at the node ${\tt t}$ in $P$.}

\item 
{\rm For each node} ${\tt t}\in\mbox{{\rm Tr}}(P)${\rm ,} $P\uarw {\tt t}$
 {\rm denotes the} subproof {\rm of} $P$ 
 {\rm whose endsequent is the sequent corresponding to the node} ${\tt t}$.

\item
{\rm 
 For each node ${\tt t}\in\mbox{{\rm Tr}}(P)$,
 $L(P\uarw{\tt t})$
 denotes the set of stage constants occurring 
in the subproof $P\uarw{\tt t}$. 

$L(P)=L(P\uarw\eps)$ denotes the set of stage constants occurring in $P$.


}
\eenu
\edf

\bdf\label{df:height}
{\rm Let} $P$ {\rm be a proof (in} $[\Pi^{0}_{1},\Pi^{0}_{1}]\mbox{{\rm -Fix}} +\fal xB(x)${\rm ) and} 
${\tt t}:\Gam$ 
{\rm a node in the proof tree ${\rm Tr}(P)$.}
{\rm We define the} \textit{height} $h({\tt t})=h({\tt t};P)\in\ome$ {\rm in} $P$ {\rm as follows:}
\benu
\item 
$h(\eps;P)=hl(\eps;P)=0$ {\rm if} $\eps:\Gam$ {\rm is the endsequent of} $P$.

{\rm In what follows let ${\tt t}:\Gam$ be an upper sequent of a rule $J$ with its lower sequent ${\tt s}:\Del$:}
\[
\infer[J]{{\tt s}:\Del}
{
\cdots
&
{\tt t}:\Gam
&
\cdots
}
\]
\item 
$h({\tt t})=h({\tt s})+1$ {\rm if $J$ is an $(h)$.}

\item 
$h({\tt t})=h({\tt s})$ {\rm otherwise.}

\eenu
\edf

In a proof $P$, each lowest rule $(h)$
with $h({\tt t};P)=0$ for its lowersequent ${\tt t}:\Gam$
 is denoted $(D)$, cf. Definition \ref{df:hydraassign}.\ref{df:hydraassign.10D}.

\bdf\label{df:labelheight}
{\rm 
Let $P$ be a proof (in} $[\Pi^{0}_{1},\Pi^{0}_{1}]\mbox{{\rm -Fix}} +\fal xB(x)${\rm ).
We define the \textit{label height} $lh({\tt t})=lh({\tt t};P)\in L\cup\{0,\mu\}$ of nodes ${\tt t}:\Gam$ 
in $P$ as follows.
\benu
\item 
$lh(\eps;P)=0$ if $\eps:\Gam$ is the endsequent of $P$.
\item
$lh({\tt t})=\mu$ if $h({\tt t})>0$.

In what follows let ${\tt t}:\Gam$ be an upper sequent of a rule $J$ with its lower sequent ${\tt s}:\Del$
such that $h({\tt t})=0$:
\[
\infer[J]{{\tt s}:\Del}
{
\cdots
&
{\tt t}:\Gam
&
\cdots
}
\]

\item
$lh({\tt t})=\max\{lh({\tt s}), lq(\vphi)\}$
if $J$ is one of basic rules, $(res)$ and $(cut)$,
where $\vphi$ denotes the minor formula of $J$ when $J$ is one of the basic rules and $(res)$, 
and $\vphi$ is the cut formula when $J$ is a $(cut)$.

\item
$lh({\tt t})=lh({\tt s})$ otherwise.
\eenu
}
\edf

\bdf\label{df:hydraassign}
{\rm 
Let $P$ be a proof in $[\Pi^{0}_{1},\Pi^{0}_{1}]\mbox{{\rm -Fix}} +\fal xB(x)$.
Let $\Ome$ be an assignment of a hydra $\Ome({\tt t})=\Ome({\tt t};P)\in \calh$ 
to each occurrence of a sequent ${\tt t}:\Gam$ in $P$. 
Also $\Ome$ assigns a label $\Ome({\tt s})\in L\cup\{0\}$
to rules ${\tt s}:(D)$.
From the assignment $\Ome$, its \textit{fixed part} $\Ome_{f}({\tt t}):=(\Ome({\tt t}))_{f}$ is determined by
Definition \ref{df:Lb}.

 If the assignment $\Ome$ enjoys the following conditions, 
then we say that $\Ome$ is a \textit{hydra assignment} for $P$.
For simplicity we write $\Ome({\tt t})$ for $\Ome({\tt t};P)$.
\benu
\item\label{df:hydraassign.1}
$\Ome({\tt t})=1$ for each axiom ${\tt t}:\Gam$.
\\

\noindent
Assume that ${\tt t}:\Gam$ is the lower sequent of a rule ${\tt s}: J$ 
and $\{{\tt t}_{i}:\Gam_{i}\}_{i<m}\, (m=1,2,3)$ denote the upper sequents of $J$. 
\[
\infer[{\tt s}:J]{{\tt t}:\Gam}
{
{\tt t}_{0}:\Gam_{0}
&
{\tt t}_{1}:\Gam_{1}
&
{\tt t}_{2}:\Gam_{2}
}
\]
\item \label{df:hydraassign.2}
$\Ome({\tt t})=\Ome({\tt t}_{0})$ {\rm if} $J$ is one of rules
$(b\fal)^{\Natural}, (\fal)^{\Natural}, (b\fal)^{\mathcal{O}}, (\fal)^{\mathcal{O}}, (\lnot I), (c)^{A}_{B}, (\lor)$ and 
$(b\exi)^{\Natural}$.

\item \label{df:hydraassign.3}
$\Ome({\tt t})=\Ome({\tt t}_{0})+\Ome({\tt t}_{1})$ if $J$ is $(\land)$.

\item \label{df:hydraassign.4}
$\Ome({\tt t})=\Ome({\tt t}_{0})+H$
 for a non-zero hydra $H\neq 0$
if $J$ is one of rules $(\exi)^{\Natural}, (b\exi)^{\mathcal{O}}, (\exi)^{\mathcal{O}}, (I)$.
In this case we write, e.g., $(I)_{H}$ for the rule $(I)$.
\\

Let $P_{\Gam,\vphi}$ be a canonically constructed proof of $\Gam,\lnot\vphi,\vphi$ using logical axioms and rules
$(\lor), (\land)$, $(b\exi)^{\Natural},(b\fal)^{\Natural}$, $(\exi)^{\Natural},(\fal)^{\Natural}$, $(b\exi)^{\mathcal{O}},(b\fal)^{\mathcal{O}}$, 
$(\exi)^{\mathcal{O}}, (\fal)^{\mathcal{O}}$
in Definition \ref{df:tautology}.
Then let $\alp_{\vphi}$ denote the (finite) ordinal canonically associated to $P_{\Gam,\vphi}$.

Namely
$H_{\vphi}=1$ if $\vphi$ is an a.p.f. or a s.p.f. or a formula of the shape $t\in I^{s}$.
$H_{\vphi}=H_{\tht_{0}}+H_{\tht_{1}}$ if $\vphi\equiv(\tht_{0}\lor\tht_{1})$.
$H_{\vphi}=H_{\tht}+4$ 
if $\vphi\equiv(\exi a<  t\, \tht(a))$.
$H_{\vphi}=H_{\tht}+1$ if $\vphi\equiv(\exi x<t\,\tht(x)), (\exi x\, \tht(x)), (\exi a\, \tht(a))$.
\\

\item \label{df:hydraassign.5}
$\Ome({\tt t})=H_{0}+\Ome({\tt t}_{0})+H_{1}$ if $J$ is a ${}_{H_{0}}(pad)_{H_{1}}$.

\item \label{df:hydraassign.5.5}
$\Ome({\tt t})=\Ome({\tt t}_{0})+H$ if $J$ is a $(res)_{H}$.

\item \label{df:hydraassign.6}
Let $J$ be a $(cut)$ with the cut formula $\tht$.

$\Ome({\tt t})=\left\{
\begin{array}{ll}
\vphi(\rk(\tht); \Ome({\tt t}_{0})+\Ome({\tt t}_{1})) & \mbox{{\rm if }} h({\tt t})=0 \spand 0\neq\rk(\tht)< \mu
\\
\Ome({\tt t}_{0})+\Ome({\tt t}_{1}) & \mbox{{\rm otherwise}}
\end{array}
\right.
$

\item \label{df:hydraassign.6.5}
Let $J$ be a $(rank)_{\mbox{\scriptsize \boldmath$C$}}$.

$\Ome({\tt t})=\left\{
\begin{array}{ll}
\vphi(\mbox{\boldmath$C$};\Ome({\tt t}_{0})) & \mbox{{\rm if }} h({\tt t})=0 
\\
\Ome({\tt t}_{0}) & \mbox{{\rm otherwise}}
\end{array}
\right.
$

\item\label{df:hydraassign.7}
Let $J$ be a $(VJ)$ with the induction term $t$.
$\Ome({\tt t})=(\Ome({\tt t}_{1})+1)\cdot mj(t)$
where $\Ome({\tt t}_{1})=\Ome({\tt t}_{0})+\Ome({\tt t}_{2})<\ome$, $mj(t)=*_{\ome}$ if $t$ is a variable.
Otherwise $t$ is a numeral $\bar{n}$. Then $mj(t)\in\{1+n,*_{\ome}\}$.

\item\label{df:hydraassign.8}
 Let $J$ be a $(TJ)$ with the induction formula $\vphi(a)$ and the induction term $s$.
$\Ome({\tt t})=(\Ome({\tt t}_{0})+\Ome({\tt t}_{1}))\cdot mj(s)$,
where $\Ome({\tt t}_{1})=H_{\vphi}$, $mj(s)=*_{\mu}$ if $s$  is a variable.
Otherwise $s$ is a constant $A$, and $mj(s)\in\{A,*_{\mu}\}$.

\[
\infer[(TJ)]{t\not<  s, \Gam}
{
\Gam,\lnot\fal b<  a\, \vphi(b), \vphi(a)
&
\lnot \vphi(t),\Gam
}
\]

\item\label{df:hydraassign.9}
$\Ome({\tt t})=\{\mu\}(\Ome({\tt t}_{0})+\Ome({\tt t}_{1}))$
if $J$ is a $(Cl.\mu)$.

\item\label{df:hydraassign.9.5}
If $J$ is a $(Cl.B)$ with $B\neq\mu$, then
$\Ome({\tt t})=\{B^{*}\}(\Ome({\tt t}_{0})+\Ome({\tt t}_{1}))$
where $B^{*}\in\{B,*_{\mu}\}$.

\item \label{df:hydraassign.10}
$\Ome({\tt t})=\ome(\Ome({\tt t}_{0}))$
if $J$ is an $(h)$ with $h({\tt t})>0$.

\item \label{df:hydraassign.10D}
$\Ome({\tt t})=D(\Ome({\tt s});\Ome({\tt t}_{0}))$ for $\Ome({\tt s})\in L^{*}$
if $J$ is a $(D)$, i.e., an $(h)$ with $h({\tt t})=0$.
In this case the rule $(D)$ is denoted by $(D_{{\tt s}})$ or by
$(D_{\mbox{\scriptsize \boldmath$C$}})$ with $\mbox{\boldmath$C$}=\Ome({\tt s})$.

\eenu

For a hydra assignment $o$ for a proof $P$ we set $\Ome(P)=\Ome(\eps:\Gam_{end})$ 
with the endsequent $\eps:\Gam_{end}$ of $P$.
}

\edf

For hydras and labels $H,H_{0}$ and labels $B\in L\cup\{0,\mu\}$,
\beqnarrs
H\ll_{B}H_{0} & :\Lrarw & o(H_{0})<o(H) \land \fal \tau\geq o(B)(K_{\tau}o(H)<d_{\tau}o(H_{0}))
\\
H\ull_{B}H_{0} & :\Lrarw & H=H_{0} \lor H\ll_{B}H_{0}
\\
H\ll_{B^{+}}H_{0} & :\Lrarw & o(H_{0})<o(H) \land \fal \tau>o(B)(K_{\tau}o(H)<d_{\tau}o(H_{0}))
\\
H\ull_{B^{+}}H_{0} & :\Lrarw & H=H_{0} \lor H\ll_{B^{+}}H_{0}
\eeqnarrs

\bdf\label{df:proofha}
{\rm 
Let $P$ be a proof in $[\Pi^{0}_{1},\Pi^{0}_{1}]\mbox{{\rm -Fix}} +\fal xB(x)$ ending with the empty sequent,
$\Ome$ a hydra assignment for $P$. 
Also $lb$ is a finite set of labels.

We say that the triple $(P,\Ome,lb)$ is a \textit{regular proof}
if the following conditions are fulfilled:}
\bdes
\item[(p0)]
{\rm 
Let ${\tt t}:\Gam$ be an uppersequent of a $(cut)$ with its cut formula $\tht$ in $P$.
Then
$\rk(\tht)<\mu+h({\tt t};P)$.

For a lower sequent ${\tt t}:\Gam$ of rules $(VJ), (TJ)$ in $P$,
$h({\tt t};P)>0$.

For a lower sequent ${\tt t}:\Gam$ of rules $(Cl.B)$ with $B^{*}\in\{\mu,*_{\mu}\}$,
$h({\tt t};P)>0$.
}

\item[(p1)]

{\rm 
Let {\tt t} be a node such that $h({\tt t};P)=0$ and {\tt u} a leaf (an axiom) in $P$ above {\tt t},
i.e., ${\tt t}\subset_{e}{\tt u}$. Then there exists an ${\tt s}:(D)$ between {\tt t} and {\tt u},
${\tt t}\subset_{e}{\tt s}\subset_{e}{\tt u}$.
In particular $\Ome(P)\in\calh_{0}$.

$L(P)\subset lb$.

For a label $C\in L$, let {\tt t} be a lowest node such that $h({\tt t};P)=0$ and $lh({\tt t};P)=C$, and $A\in L(P\uarw{\tt t})$
be a label occurring above {\tt t} such that $C\leq A$.
Then $A\in\Ome_{f}({\tt t})$.
}

\item[(p2)]
{\rm 
Let $\infer[{\tt s}:(c)^{B}_{A}]{\Gam_{0}}{{\tt t}:\Gam_{1};H}$ be a rule in $P$.

 
 

Then
$L(P\uarw{\tt s})\ll_{B}H_{0}$ and  $H\ull_{B}H_{0}$
for $H=\Ome({\tt t};P)$ and $A=d_{B}(H_{0})$.
 }
 
\item [(p3)] 
{\rm 
Let $\infer[{\tt s}:(D)]{{\tt t}_{1}:\Gam}{{\tt t}_{0}:\Gam}$ be a rule in $P$. Then
$L(P\uarw {\tt s})\ull_{\mu} \Ome({\tt s})$.
}

\edes
\edf
Note that by {\bf (p2)} $L(P\uarw{\tt s})\cap B<C$ holds for rules ${\tt s}:(c)^{B}_{A}$, 
i.e., any constant $C$ occurring above the rule ${\tt s}$ is $C< A$ if $C< B$

 Observe again that the relation 
 \beqnarrs
 && `x \mbox{ is a triple } (p,\Ome,lb) \mbox{ such that } p \mbox{ is a proof in the system } [\Pi^{0}_{1},\Pi^{0}_{1}]\mbox{{\rm -Fix}} +\fal xB(x) \\
 && \mbox{ with an h.a } \Ome \spand \Ome(p)=H"
 \eeqnarrs
 is elementary recursive.

\bprp \label{prp:4}
Assume $[\Pi^{0}_{1},\Pi^{0}_{1}]\mbox{{\rm -Fix}} +\fal xB(x)$ is inconsistent. 
Then there exists a regular proof $(P,\Ome,\emptyset)$.
\eprp
\bprf
Let $P_{0}$ be a proof in $[\Pi^{0}_{1},\Pi^{0}_{1}]\mbox{{\rm -Fix}} +\fal xB(x)$ ending with the empty sequent. 
Leaves for complete induction schema are replaced by the following:
{\scriptsize
\[
\infer[(\fal),(\lor)]{\Gam,A(0)\land \fal y(A(y)\to A(y^{\prime}))\to\fal x\, A(x); (2d+2)\cdot*_{\ome}}
{
  \infer[(VJ)]{\Gam,\Del,A(x); (2d+2)\cdot*_{\ome}}
   {
    \infer*{\Gam,\lnot A(0),A(0);d}{}
    &
    \infer[(\land),(\exi)]{\Gam,\Del,\lnot A(y), A(y^{\prime}); 2d+1}
    {
    \infer*{\Gam,A(y),\lnot A(y);d}{}
    &
    \infer*{\Gam,\lnot A(y^{\prime}),A(y^{\prime}) ; d}{}
    }
    & 
    \infer[(pad)_{1}]{\Gam,\lnot A(x),A(x); d+1}
    {\infer*{\Gam,\lnot A(x),A(x);d}{}}
    }
}
\]
}
where $\Del=\{\lnot A(0),\exi y(A(y)\land A(y^{\prime}))\}$,
$d=H_{A(y)}$, and $*_{\ome}=mj(x)$.

Leaves for transfinite induction schema are replaced by
{\scriptsize
\[
\hspace{-10mm}
\infer[(\fal),(\lor)]{\Gam,\fal b(\fal a< b\, A(a)\to A(b))\to\fal b\, A(b); (3d+5)\cdot*_{\mu}+2d+5}
{
\infer[(cut)]{\Gam,\lnot Prg, A(b); d_{1}\cdot*_{\mu}+d_{0}}
{
 \infer[(b\fal)]{\Gam,\lnot Prg,\fal a< b\, A(a)}
  {
  \infer[(TJ)]{a\not< b, \Gam,\Del; d_{1}\cdot*_{\mu}}
   {
    \infer[(\land),(\exi)]{\Gam,\lnot Prg,\lnot\fal a< b A(a), A(b); d_{0}}
    {
    \infer*{\Gam,\fal a< b\, A(a),\lnot\fal a< b A(a);d+4}{}
    &
    \infer*{\Gam,\lnot A(b),A(b) ; d}{}
    }
    & 
    \hspace{-2mm}
   \infer*{\Gam,\Del,A(a),\lnot A(a);d}{}
    }
   }
  &
  \hspace{-25mm}
  \infer*{\Gam,\lnot Prg,\lnot\fal a< b\, A(a),A(b); d_{0}}{}
 }
}
\]
}
where $\Del=\{\lnot Prg, A(a)\}$ with $Prg\equiv(\fal b(\fal a< b\, A(a)\to A(b)))$ and 
$d=H_{A(a)}$, $d+4= H_{\fal a< b\, A(a)}$, $d_{0}=2d+5$, and
$d_{1}=3d+5$.
Also $*_{\mu}=mj(b)$.

$P_{0}$ contains none of rules $(c),(pad),(h)$, and no constant of stage sort occurs in $P_{0}$ besides $0^{\mathcal{O}}$. 
Below the endsequent of $P_{0}$ attach some $(h)$'s to enjoy the condition {\bf (p1)}.
A hydra assignment $\Ome$ for $P$ is chosen canonically, and $\Ome({\tt s})=\emptyset$.
Namely the bottom of $P$ looks like
\[
P=
\left.
\begin{array}{c}
\infer[{\tt s}:(D)]{\eps:\emptyset;D_{\emptyset}(H)}
{
 \infer*{{\tt t}:\emptyset;H}
 {
  \infer[(h)]{\emptyset}
 {
 \infer*[P_{0}]{\emptyset}{}
 }
 }
}
\end{array}
\right.
\]
The resulting quadruple $(P,\Ome,\emptyset)$ is regular.
\eprf

\subsection{Inversions}\label{subsec:inversion}

Let $\Ome$ be a hydra assignment for a proof $P$, and ${\tt t}:\Gam,\tht_{1}$ a node in $P$.
\[
P=
\left.
\begin{array}{c}
\infer*{}
{
 \infer*[P_{0}]{{\tt t}:\Gam,\tht_{1}:H}{}
 }
\end{array}
\right.
\]
Let us define a proof $P_{0}^{\prime}$ of a $\Gam,\tht_{0}$ by inversion so that
$\Ome({\tt t};P_{0}^{\prime})=H=\Ome({\tt t};P)$ according to the formulas $\tht_{1}$.
\benu
\item
$\tht_{1}\equiv(\fal x\leq \bar{n}\, \tht(x))$ for the bounded number quantifier $\fal x\leq\bar{n}$:
For $k\leq n$, let $\tht_{0}\equiv\tht(\bar{k})$.
To get a $P_{0}^{\prime}$ by inversion,
change $(b\fal)^{\Natural}$ to $(pad)_{0}$, and eliminate the false $\bar{k}\not\leq\bar{n}$ if necessary:
\[
\infer[(b\fal)^{\Natural}]{\Psi,\fal x\leq \bar{n}\, \tht(x); H_{0}}
{
\Psi,\fal x\leq\bar{n}\tht(x),y\not\leq\bar{n}\,\tht(y);H_{0}
}
\leadsto
\infer[(pad)_{0}]{\Psi,\tht(\bar{k}); H_{0}}
{
\Psi,\tht(\bar{k});H_{0}
}
\]
and
\[
\infer[(VJ)]{\Phi; (H_{2}+1)\cdot *_{\ome}}
 {
   \infer*{\Phi,\vphi(0); H_{1}}{}
   &
   \infer*{\lnot \vphi(x),\Phi, \vphi(x') ; H_{2}}{}
   &
   \infer*{\lnot\vphi(y),\Phi; H_{3}}{}
  }
\]
turns to the following with $mj(\bar{k})=*_{\ome}$.
\[
\infer[(VJ)]{\Phi; (H_{2}+1)\cdot *_{\ome}}
 {
   \infer*{\Phi,\vphi(0); H_{1}}{}
   &
   \infer*{\lnot \vphi(x),\Phi, \vphi(x') ; H_{2}}{}
   &
   \infer*{\lnot\vphi(\bar{k}),\Phi; H_{3}}{}
  }
\]
The resolvent $\Pi_{m}=\Pi^{\prime}_{m}\cup\{\fal x\leq\bar{n}\,\tht(x)\}$ of a $\vphi$ turns to a resolvent
$\Pi^{\prime}_{m}\cup\{\tht(\bar{k})\}$ of the same formula.
{\small
\[
\infer[(res)_{K}]{\bigcup_{i\leq p}\Pi_{i},\Del_{0}}{t\leq p\land \vphi(t),\bigcup_{i\leq p}\Pi_{i},\Del_{0}}
\leadsto
\infer[(res)_{K}]{\bigcup_{i\neq m}\Pi_{i},\Pi^{\prime}_{m}\cup\{\tht(\bar{k})\},\Del_{0}}{t\leq p\land \vphi(t),\bigcup_{i\neq m}\Pi_{i},\Pi^{\prime}_{m}\cup\{\tht(\bar{k})\},\Del_{0}}
\]
}
Moreover when the variable $y$ in a $(res)$ is replaced by $\bar{k}\leq\bar{m}$, 
one of the formulas $\vphi(\bar{k})$ and $\vphi(\ovl{k+1})$ is replaced by its resolvent $\Pi_{k}$ and $\Pi_{k+1}$ by inversions, resp.
\[
 \infer[(res)_{K}]{\bigcup_{k\leq m}\Pi_{k},\Del_{0};H_{0}+K}
{
 y\leq \bar{m}\land\vphi(y),\bigcup_{k\leq m}\Pi_{k},\Del_{0};H_{0}
 }
\leadsto
\infer[(pad)_{K}]{\bigcup_{k\leq m}\Pi_{k},\Del_{0};H_{0}+K}
{
\bigcup_{k\leq m}\Pi_{k},\Del_{0};H_{0}
 }
\]

\item
$\tht_{1}\equiv(\fal x\tht(x))$ for the unbounded number quantifier $\fal x$:
Similar to the case for bounded universal number quantifiers, but there is no concern with resolvents.

\item
$\tht_{1}\equiv(\fal a\tht(a))$ for the stage quantifier $\fal a$:
For a stage constant $C$, let $\tht_{0}\equiv\tht(C)$.
To get a $P_{0}^{\prime}$ by inversion,
change $(\fal)^{\mathcal{O}}$ to $(pad)_{0}$ if necessary:
\[
\infer[(\fal)^{\mathcal{O}}]{\Psi,\fal a\tht(a); H_{0}}
{
\Psi,\fal a \tht(a),\tht(a_{0});H_{0}
}
\leadsto
\infer[(pad)_{0}]{\Psi,\tht(A); \bet_{0}}
{
\Psi,\tht(A);\bet_{0}
}
\]
and
\[
 \infer[(TJ)]{s\not< a_{0},\Phi; H_{1}+(H_{\vphi}+H_{2})\cdot *_{\mu}}
 {
  \infer*{\Phi, \lnot \vphi(s);H_{1}}{}
  &
  \infer*{\vphi(c),\lnot\fal b< c\, \vphi(b),\Phi;H_{2}}{}
 }
\]
turns to the following with $mj(A)=*_{\mu}$.
\[
 \infer[(TJ)]{s\not< A,\Phi; H_{1}+(H_{\vphi}+H_{2})\cdot *_{\mu}}
 {
  \infer*{\Phi, \lnot \vphi(s);H_{1}}{}
  &
  \infer*{\vphi(c),\lnot\fal b\prec c\, \vphi(b),\Phi;H_{2}}{}
 }
\]

\item
$\tht_{1}\equiv(\fal a<A\tht(a))$ for the stage bounded quantifier $\fal a$:
Similar to the case for unbounded universal stage quantifiers.

\item
$\tht_{1}\equiv(\tht_{2}\land\tht_{3})$:
For $i=2,3$, let $\tht_{0}\equiv\tht_{i}$.
To get a $P_{0}^{\prime}$ by inversion,
change $(\land)$ to $(pad)$ if necessary:
\[
\infer[(\land)]{\Psi,\tht_{2}\land\tht_{3}; H_{2}+H_{3}}
{
\Psi,\tht_{2}\land\tht_{3},\tht_{2};H_{2}
&
\Psi,\tht_{2}\land\tht_{3},\tht_{3};H_{3}
}
\leadsto
\infer[(pad)]{\Psi,\tht_{2}; H_{2}+H_{3}}
{
\Psi,\tht_{i};H_{i}
}
\]
Moreover 
the resolvent $\Pi_{m}=\Pi^{\prime}_{m}\cup\{\tht_{2}\land\tht_{3}\}$ of a $\vphi$ turns to a resolvent
$\Pi^{\prime}_{m}\cup\{\tht_{i}\}$ of the same formula.
{\small
\[
\infer[(res)_{K}]{\bigcup_{i\leq p}\Pi_{i},\Del_{0}}{t\leq \bar{p}\land \vphi(t),\bigcup_{i\leq p}\Pi_{i},\Del_{0}}
\leadsto
\infer[(res)_{K}]{\bigcup_{i\neq m}\Pi_{i},\Pi^{\prime}_{m}\cup\{\tht_{i}\},\Del_{0}}{t\leq \bar{p}\land \vphi(t),\bigcup_{i\neq m}\Pi_{i},\Pi^{\prime}_{m}\cup\{\tht_{i}\},\Del_{0}}
\]
}

\item
$\tht_{1}\equiv(\exi x\leq \bar{n}\tht(x))$:
Let $\tht_{0}=\{\tht(\bar{k}): k\leq n\}$.
To get a $P_{0}^{\prime}$ by inversion,
change $(b\exi)^{\Natural}$ to $(pad)$ if necessary: for $k\leq n$,
{\small
\[
\infer[(b\exi)^{\Natural}_{K}]{\exi x\leq \bar{n}\tht(x),\Del_{0};H_{0}+K}
{
 \bar{k}\leq \bar{n}\land\tht(\bar{k}),\exi x\leq \bar{n}\tht(x),\Del_{0};H_{0}
 }
\leadsto
\infer[(pad)_{K}]{\{\tht(\bar{k}):k\leq n\},\Del_{0};H_{0}+K}
{
 \{\tht(\bar{k}):k\leq n\},\Del_{0};H_{0}
 }
\]
}
where
some rules $(\land)$ with the principal formula $\bar{k}\leq\bar{n}\land\tht(\bar{k})$
is also replaced by paddings together with eliminating the left upper part of the $(\land)$ if $k\leq n$,
and eliminating the left upper part of the $(\land)$ if $k>n$.

If the witnessing term $t$ is a variable $y$, then the rule becomes a $(res)$:
{\small
\[
\infer[(b\exi)^{\Natural}_{K}]{\exi x\leq \bar{n}\tht(x),\Del_{0};H_{0}+K}
{
y\leq \bar{n}\land\tht(y),\exi x\leq \bar{n}\tht(x),\Del_{0};H_{0}
 }
\leadsto
\infer[(res)_{K}]{\{\tht(\bar{k}):k\leq n\},\Del_{0};H_{0}+K}
{
 y\leq \bar{n}\land\tht(y),\{\tht(\bar{k}):k\leq n\},\Del_{0};H_{0}
 }
\]
}
The case when $t\equiv y^{\prime}$ is similar:
{\small
\[
\infer[(b\exi)^{\Natural}_{K}]{\exi x\leq \bar{n}\tht(x),\Del_{0};H_{0}+K}
{
y^{\prime}\leq \bar{n}\land\tht(y^{\prime}),\exi x\leq \bar{n}\tht(x),\Del_{0};H_{0}
 }
\leadsto
\infer[(res)_{K}]{\{\tht(\bar{k}):k\leq n\},\Del_{0};H_{0}+K}
{
 y^{\prime}\leq \bar{n}\land\tht(y^{\prime}),\{\tht(\bar{k}):k\leq n\},\Del_{0};H_{0}
 }
\]
}

Moreover
the resolvent $\Pi_{m}=\Pi^{\prime}_{m}\cup\{\exi x\leq\bar{n}\tht(x)\}$ of a $\vphi$ turns to a resolvent
$\Pi^{\prime}_{m}\cup\{\tht(\bar{k}):k\leq n\}$ of the same formula.
{\small
\[
\infer[(res)_{K}]{\bigcup_{i\leq p}\Pi_{i},\Del_{0}}{t\leq \bar{p}\land \vphi(t),\bigcup_{i\leq p}\Pi_{i},\Del_{0}}
\leadsto
\infer[(res)_{K}]{\bigcup_{i\neq m}\Pi_{i},\Pi^{\prime}_{m}\cup\{\tht(\bar{k}):k\leq n\},\Del_{0}}{t\leq \bar{p}\land \vphi(t),\bigcup_{i\neq m}\Pi_{i},\Pi^{\prime}_{m}\cup\{\tht(\bar{k})\},\Del_{0}}
\]
}

\item
$\tht_{1}\equiv(\tht_{2}\lor\tht_{3})$:
Let $\tht_{0}=\{\tht_{2},\tht_{3}\}$.
To get a $P_{0}^{\prime}$ by inversion,
change $(\lor)$ to $(pad)$ if necessary: 
\[
\infer[(\lor)_{K}]{\tht_{2}\lor\tht_{3},\Del_{0};H_{0}+K}
{
 \tht_{2}\lor\tht_{3},\tht_{i},\Del_{0};H_{0}
 }
\leadsto
\infer[(pad)_{K}]{\tht_{2},\tht_{3},\Del_{0};H_{0}+K}
{
\tht_{2},\tht_{3},\Del_{0};H_{0}
 }
\]

Moreover
the resolvent $\Pi_{m}=\Pi^{\prime}_{m}\cup\{\tht_{2}\lor\tht_{3},\}$ of a $\vphi$ turns to a resolvent
$\Pi^{\prime}_{m}\cup\{\tht_{2},\tht_{3}\}$ of the same formula.
{\small
\[
\infer[(res)_{K}]{\bigcup_{i\leq p}\Pi_{i},\Del_{0}}{t\leq \bar{p}\land \vphi(t),\bigcup_{i\leq p}\Pi_{i},\Del_{0}}
\leadsto
\infer[(res)_{K}]{\bigcup_{i\neq m}\Pi_{i},\Pi^{\prime}_{m}\cup\{\tht_{2},\tht_{3}\},\Del_{0}}{t\leq \bar{p}\land \vphi(t),\bigcup_{i\neq m}\Pi_{i},\Pi^{\prime}_{m}\cup\{\tht(\bar{k})\},\Del_{0}}
\]
}
\eenu

\section{Rewritings}\label{sec:rewriting}
In this section we define rewritings on proofs in such a way that each rewriting corresponds to 
a move on hydras attached to proofs.

Let $P$ be a proof in $[\Pi^{0}_{1},\Pi^{0}_{1}]\mbox{{\rm -Fix}} +\fal xB(x)$.
$num(P)$ denotes the set of numerals $\bar{n}$ occurring in $P$, and
$Fml(P)$ denotes the set of formulas occurring in $P$.
For formulas $\vphi$, let $q(\vphi)$ denote the number of occurrences of logical connectives 
$\land,\lor,\fal,\exi$ in $\vphi$.
Then let
$c(P)=\max(\{n+1: \bar{n}\in num(P)\}\cup\{q(\vphi): \vphi\in Fml(P)\}\cup\{d_{\cala}+1\})$.

Note that for $\vphi\in Fml(P)$, there exists a label $A\in lb(\vphi)\cup\{0,\mu\}$ such that
$\rk(\vphi)< A+c(P)$ by Lemma \ref{lem:rank}.\ref{lem:rank15}.

Next for terms $t$ of number sort,
$c^{\prime}(t)$ denotes a natural number defined as follows.
$c^{\prime}(\bar{n})=0$ for numerals $\bar{n}$ and $c^{\prime}(t^{\prime})=c^{\prime}(t)+1$
if $t$ is not a numeral.
Let $w(P)$ denote the set of witnessing terms of rules $(b\exi)^{\Natural},(\exi)^{\Natural}$
and induction terms of rules $(VJ)$ in $P$.
Then let $c^{\prime}(P)=\max(\{0\}\cup\{c^{\prime}(t): t \in w(P)\})$.

Suppose $c^{\prime}(P)\leq 1$, and let $P_{0}$ be a proof obtained from $P$ by substituting
a numeral $\bar{n}$ for a variable of number sort with $n<c(P)$.
We see then that $c^{\prime}(P_{0})\leq c^{\prime}(P)$ and $c(P_{0})\leq c(P)+1$.

Let $(P_{-2},\Ome_{-2},\emptyset)$ be a regular proof in $[\Pi^{0}_{1},\Pi^{0}_{1}]\mbox{{\rm -Fix}} +\fal xB(x)$ 
in which no stage constant except $0^{\mathcal{O}}$ occurs, cf.\,Proposition \ref{prp:4}.
Suppose $c^{\prime}(P_{-2})>1$.
Let us construct a regular proof $P_{-1}$ without stage constant except $0^{\mathcal{O}}$
such that $c^{\prime}(P_{-1})\leq 1$.
Such a proof $P_{-1}$ is obtained by replacing rules $(b\exi)^{\Natural},(\exi)^{\Natural},(VJ)$ with
a `big' term $t^{\prime}$ with $c^{\prime}(t)> 1$ repeatedly as follows.
Replace
\[
\infer[(\exi)^{\Natural}]{\exi x\vphi(x),\Gam}{\exi x\vphi(x),\vphi(t^{\prime}),\Gam}
\]
with a fresh variable $y$ by the following:
\[
\infer[(cut)]{\exi x\vphi(x),\Gam}
{
\infer[(\exi)^{\Natural}]{\exi y(y=t)}{t=t}
&
\infer[(\fal)^{\Natural}]{\exi x\vphi(x),\Gam,\fal y(y\neq t)}
{
\infer[(\exi)^{\Natural}]{\exi x\vphi(x),\Gam,y\neq t}
{
 \infer[(cut)]{\exi x\vphi(x),\vphi(y^{\prime}),\Gam,y\neq t}
 {
  \exi x\vphi(x),\vphi(t^{\prime}),\Gam
  &
  \infer[(cut)]{\lnot\vphi(t^{\prime}),\vphi(y^{\prime}),y\neq t}
  {
   \infer*{\lnot\vphi(t^{\prime}),\vphi(y^{\prime}),y^{\prime}\neq t^{\prime}}{}
   &
   y^{\prime}=t^{\prime},y\neq t
   }
  }
 }
}
}
\]
Replace
\[
\infer[(VJ)]{\Gam}{\Gam, \vphi(0) & \lnot \vphi(x),\Gam, \vphi(x') & \lnot \vphi(t^{\prime}),\Gam}
\]
by the following with a fresh variable $y$:
\[
\infer[(cut)]{\Gam}
{
\infer[(\exi)^{\Natural}]{\exi y(y=t)}{t=t}
&
\hspace{-15mm}
\infer[(\fal)^{\Natural}]{\Gam,\fal y(y\neq t)}
{
\infer[(VJ)]{\Gam,y\neq t}
{
\Gam, \vphi(0) 
& 
\lnot \vphi(x),\Gam, \vphi(x') 
& 
\hspace{-10mm}
\infer[(cut)]{\lnot \vphi(y^{\prime}),\Gam,y\neq t}
{
 \lnot \vphi(t^{\prime}),\Gam
 &
   \infer[(cut)]{\vphi(t^{\prime}),\lnot\vphi(y^{\prime}),y\neq t}
  {
   \infer*{\vphi(t^{\prime}),\lnot\vphi(y^{\prime}),y^{\prime}\neq t^{\prime}}{}
   &
   y^{\prime}=t^{\prime},y\neq t
   }
   }
  }
}
}
\]
Then the resulting proof $P_{-1}$ such that $c^{\prime}(P_{-1})\leq 1$
can be assumed to be regular for some $\Ome_{-1}$. 
Otherwise insert some rules $(h)$ for newly arising $(cut)$'s.
Let $c=c(P_{-1})$, and for $k\leq c$, $P_{i}$ be a proof obtained from $P_{-1}$ by adding a rule $(pad)_{c-k}$ 
as the last rule:
\[
P_{k}=
\left.
\begin{array}{c}
\infer[(pad)_{c-k}]{\emptyset;D_{\emptyset}(H)+c-k}
{
 \infer*[P_{-1}]{\emptyset;D_{\emptyset}(H)}{}
 }
 \end{array}
 \right.
 \]
$(P_{k},\Ome_{k},\emptyset)$ with $\Ome_{k}=\Ome_{-1}$ is a regular proof such that
$c^{\prime}(P_{k})\leq 1$ and $c(P_{k})\leq c$.
Moreover $(D_{\emptyset}(H)+c-k,\emptyset)\to_{k} (D_{\emptyset}(H)+c-k-1,\emptyset)$ for $0\leq k< c$.
This yields a regular proof $(P_{c},\Ome_{c},\emptyset)$ such that $c^{\prime}(P_{c})\leq 1$ and
$c(P_{c})\leq c$.

We construct regular proofs $(P[\ell],\Ome[\ell],lb[\ell])$ for $c\leq\ell<\ome$ in such a way that
$(P[c],\Ome[c],lb[c])=(P_{c},\Ome_{c},\emptyset)$, $L(P[\ell])\subset lb[\ell]$, and
$(P[\ell],lb[\ell])\to_{\ell}(P[\ell+1],lb[\ell+1])$ for each $\ell\geq c$.
Moreover  for each $\ell\geq c$
\beqn\label{eq:cprime}
c^{\prime}(P[\ell])\leq 1 \spand c(P[\ell])\leq \ell
\eeqn

This means that $\{t_{\ell}\}_{\ell}$ is an infinite path through the tree $Tr(H_{0},lb_{0})$,
where $H_{0}=\Ome_{0}(P_{0})$, $lb_{0}=\emptyset$, $t_{0}=\eps$, and 
$t_{\ell+1}=t_{\ell}*(n_{\ell})$ such that $(\Ome[\ell+1])(P[\ell]+1)$ is the $n_{\ell}$'s move from
$(\Ome[\ell])(P[\ell])$.

Let $P=P[\ell], lb=lb[\ell], \Ome=\Ome[\ell]$ and
$P^{\prime}=P[\ell+1], lb^{\prime}=lb[\ell+1], \Ome^{\prime}=\Ome[\ell+1]$.
Except {\bf Case 3.4} and {\bf Case 3.5} in subsection \ref{subsec:indrfl}
$lb^{\prime}=lb[\ell+1]=lb[\ell]=lb$ holds.


\bdf Main branch\\
{\rm 
Let $P$ be a proof ending with the empty sequent. 
The \textit{main branch} of $P$ is a series
$\{{\tt t}_{i}:\Gam_i\}_{i\leq n}$ of occurrences of sequents in $P$ such that:
\benu
\item 
${\tt t}_{0}:\Gam_{0}$ is the endsequent of $P$, i.e., ${\tt t}_{0}=\eps$.

\item 
For each $i<n$ ${\tt t}_{i+1}:\Gam_{i+1}$ {\rm is the} \textit{right} {\rm upper sequent of a rule} $J_{i}$ {\rm so that} 
${\tt t}_{i}:\Gam_{i}$ {\rm is the lower sequent of} $J_{i}$ {\rm and} $J_{i}$ {\rm is one of the rules} ${}_{H^{\prime}}(pad)_{0}, (res)_{0},(h), (c)$ and $(cut)$ with a cut formula in one of the shapes
$\exi x\vphi, \exi a< s\vphi, \exi a\vphi, t\in I^{s}$.

\item 
Either ${\tt t}_{n}:\Gam_{n}$ is an axiom, or
 ${\tt t}_{n}:\Gam_{n}$ is the lower sequent of one of the rules
$(b\exi)^{\Natural},(\exi)^{\Natural},(b\exi)^{\mathcal{O}},(\exi)^{\mathcal{O}}, (VJ), (TJ), (I), (Cl)$ 
and ${}_{H^{\prime}}(pad)_{H}$ and $(res)_{H}$
with $H\neq 0$, or
${\tt t}_{n}:\Gam_{n}$ is the lower sequent of a $(cut)$ with an \textit{unbounded}
cut formula in one of the shapes
$\vphi_{0}\lor\vphi_{1}$ or $\exi x\leq n \vphi$ for numerals $n$.
\eenu
The sequent ${\tt t}_{n}:\Gam_{n}$ is said to be the \textit{top} (of the main branch) of the proof $P$.
}
\edf

Let $\Phi$ denote the top of the proof $P$ with the hydra assignment $\Ome$. 
Observe that we can assume $\Phi$ contains no free variable
for otherwise substitute $0^{\Natural}$ for number variables, and $0^{\mathcal{O}}$ for stage variables.
The same hydra assignment works for the substituted proof.

In each case below the new hydra assignment $\Ome'$ for the new proof $P'$ is defined obviously from the hydra
assignment $\Ome$ 
and the subscripts $H$ of the displayed padding rules ${}_{H^{\prime}}(pad)_{H}$.

\subsection{Rewritings by necrosis}\label{subsec:toppad}

In this subsection we consider the cases when the top 
$\Phi$ is either  
the lower sequent of a padding $(p)_{H}={}_{H^{\prime}}(pad)_{H}$ with $H\neq 0$
or the lower sequent of one of rules one of rules 
$(p)_{H}=(b\exi)^{\Natural}_{H}, (\exi)^{\Natural}_{H}, (b\exi)^{\mathcal{O}}_{H}, (\exi)^{\mathcal{O}}_{H}, (I)_{H}$
with $H\not\in\{0,1\}$
or an axiom $(ax)$.
\\

\noindent
{\bf Case 1}. 
$\Phi$ is 
the lower sequent of a padding $(p)_{H}={}_{H^{\prime}}(pad)_{H}$ with $H\neq 0$.
Then kill the padding by (Necrosis) $(H,\emptyset)\to_{\ell}(0,\emptyset)$ in Definition \ref{df:move}.\ref{move:1}.
\[
P=
\left.
\begin{array}{c}
\infer*{}
{
 \infer[(p)_{H}]{\Phi;\, H_{0}^{\prime}+H}
 {
  \infer*{\cdots;\, H_{0}}{}
 }
} 
\end{array}
\right.
\msfiv
P':=
\left.
\begin{array}{c}
\infer*{}
{
 \infer[{(p)_{0}}]{\Phi;\, H^{\prime}_{0}}
 {
  \infer*{\cdots;\, H_{0}}{}
 }
} 
\end{array}
\right.
\]
 {\bf Case 2}.
$\Phi$ is the lower sequent of one of rules one of rules 
$(p)_{H}=(b\exi)^{\Natural}_{H}$, $(\exi)^{\Natural}_{H}$, $(b\exi)^{\mathcal{O}}_{H}$,
$(\exi)^{\mathcal{O}}_{H}$, $(cut)_{H}$
with $H\not\in\{0,1\}$.
Again kill the padding by (Necrosis) $(H,\emptyset)\to_{\ell}(1,\emptyset)$
in Definition \ref{df:move}.\ref{move:1}.
\[
P=
\left.
\begin{array}{c}
\infer*{}
{
 \infer[(p)_{H}]{\Phi;\, H_{0}+H}
 {
  \infer*{\cdots;\, H_{0}}{}
 }
} 
\end{array}
\right.
\msfiv
P':=
\left.
\begin{array}{c}
\infer*{}
{
 \infer[{(p)_{1}}]{\Phi;\, H_{0}+1}
 {
  \infer*{\cdots;\, H_{0}}{}
 }
} 
\end{array}
\right.
\]
 \\
 {\bf Case 3}.
$\Phi$ is a nonlogical axiom:
Then, since $\Phi$ contains no free variable and $\fal x B(x)$ is assumed to be true, 
$\Phi$ contains a $\tht$ which is either a true a.p.f. or a true s.p.f. or $n\not\in I^{0^{\mathcal{O}}}$.
Let $\Phi=\tht,\Del_{0}$ with $\rk(\tht)=lq(\tht)=0$.
Eliminate the false prime formula $\lnot \tht$ and insert a $(pad)_{0}$. 
$\Ome(P')$ is obtained from $\Ome(P)$ by (Necrosis), $(H+K,\emptyset)\to_{\ell}(H,\emptyset)$.
Note that $lh({\tt t};P^{\prime})=lh({\tt t};P)$ since $lq(\tht)=0$.
\[
P= 
\left.
\begin{array}{c}  
\infer*{}
{
 \infer{\Gam,\Del;\, H+K}
 {
  \infer*{{\tt t}:\Gam,\lnot \tht;\, H}{}
 &
  \infer*{\tht,\Del;\, K}{\tht,\Del_0;\, 1}
 }
}
\end{array}
\right.
\msfiv
P':=
\left.
\begin{array}{c}   
\infer*{}
{
  \infer[(pad)_{0}]{\Gam,\Del;\, H}
  {
   \infer*{{\tt t}:\Gam;\, H}{}
   }
  }
\end{array}
\right.
\]
{\bf Case 4}. 
$\Phi$ is a logical axiom: 
$\Phi=\lnot\tht,\tht,\Del_0$, where $\tht$ is an a.p.f. or a s.p.f. 
Note that the case when $\tht\equiv(n\in I^{A})$ is excluded since the endsequent is empty, and
$n\not\in I^{A}$ is not an $\exi$-formula.
Consider a $(cut)$ whose right upper sequent is a sequent $\lnot\tht,\Del$ with $\tht\in\Del$, and
$\tht$ is its cut formula. 
$\Ome(P')$ is obtained from $\Ome(P)$ by (Necrosis).
\[
P=
\left.
\begin{array}{c}    
\infer*{}
{
 \infer{\Gam,\Del;\, H+K}
 {
  \infer*{\Gam,\tht;\, H}{}
 &
  \infer*{\lnot\tht,\Del;\, K}{\lnot\tht,\tht,\Del_0; 1}
 }
}
\end{array}
\right.
\msfiv
P':=
\left.
\begin{array}{c}   
\infer*{}
{
  \infer[(pad)_{0}]{\Gam,\Del;\,H}
  {
   \infer*{\Gam,\tht;\, H}{}
   }
}
\end{array}
\right.
\]

\subsection{Rewritings on bounded logical rules}\label{subsec:topbnd}

In this subsection we consider the cases when the top $\Phi$ is a lower sequent of a $(cut)$ with an
\textit{unbounded} cut formula 
in one of the shapes
$\vphi_{0}\lor\vphi_{1}$ or $\exi x\leq \bar{n} \vphi$ for numerals $\bar{n}$.
Let us consider the latter case, and $P$ be the following.
\[
P=
\left.
\begin{array}{c}
\infer*{}
{
\infer[{\tt s}]{{\tt u}:\Gam_{3};\ome(H_{3}+1)}
{
\infer*{\Gam_{3};H_{3}+1}
{
\infer{{\tt t}:\Gam,\Del;H_{2}}
 {
  \infer*[P_{0}]{\Gam,\fal x\leq \bar{n}\lnot\vphi(x);H_{0}}{}
  &
  \infer*[P_{1}]{\exi x\leq \bar{n}\vphi(x),\Del;H_{1}}{}
 }
} 
}
}
\end{array}
\right.
\]
where
$H_{2}=H_{0}+H_{1}+1$, and
{\tt s} denotes the uppermost rule $(h)$ below the top {\tt t}.
We have $h({\tt u})>0$ by {\bf (p0)} since $\vphi$ is unbounded.
We obtain $(\ome(H_{3}+1),lb)\to_{\ell}(\ome(H_{3})\cdot n,lb)$ by Definition \ref{df:move}.\ref{move:d},
where $n<\ell$ by (\ref{eq:cprime}).

For each $k\leq n$,
let $P_{0}^{\prime}(k)$ be a proof of $\Gam,\lnot\vphi(k)$, which is obtained from $P_{0}$ by inversion.
Let $P_{1}^{\prime}$ be obtained from $P_{1}$ by
replacing the formula $\exi x\leq n\vphi(x)$ by the set $\{\vphi(k):k\leq n\}$, cf.\,subsection \ref{subsec:inversion}.
Let
$P^{\prime}$ be the following when $\vphi$ is an unbounded $\exi$-formula.
\[
\infer*{}
{
\infer[J]{\Gam_{3};\ome(H_{3})(n+2)}
{
\infer{\Gam_{3},\lnot\vphi(0);\ome(H_{3})}
{
\infer*{\Gam_{3},\lnot\vphi(0);H_{3}}
{
  \infer{\Gam,\Del,\lnot\vphi(0);H_{0}+H_{1}}
  {
  \infer*[P_{0}^{\prime}(0)]{\Gam,\lnot\vphi(0);H_{0}}{}
  }
  }
 }
&
\cdots
&
\infer{\Gam_{3},\lnot\vphi(n);\ome(H_{3})}
{
 \infer*{\Gam_{3},\lnot\vphi(n);H_{3}}
{
  \infer{\Gam,\Del,\lnot\vphi(n);H_{0}+H_{1}}
  {
  \infer*[P_{0}^{\prime}(n)]{\Gam,\lnot\vphi(n);H_{0}}{}
  }
 }
 }
&
\infer{\{\vphi(k):k\leq n)\},\Gam_{3};\ome(H_{3})}
{
 \infer*{\{\vphi(k):k\leq n\},\Gam_{3};H_{3}}
 {
 \infer{\{\vphi(k):k\leq n\},\Gam,\Del;H_{0}+H_{1}}
  {
   \infer*[P_{1}^{\prime}]{\{\vphi(k):k\leq n\},\Del;H_{1}}{}
  }
 } 
 }
}
}
\]
where $J$ denotes several $(cut)$'s with unbounded cut formulas $\vphi(k)$.
\\

\noindent
Note that due to the the rewritings in this subsection, we can assume the following.
Let $\vphi$ is an unbounded formula such that $\vphi^{B}$ is in the upper sequent of a rule $(c)^{B}_{A}$ and
$\vphi^{A}$ is in its lower sequent.
Then $\vphi^{B}$ is one of the formulas
$\fal y \calb_{i}(I^{< B},\bar{n},y)$, $\Sig^{B}$-formula $\bar{m}\in I^{< B}$, or
$\Pi^{B}$-formula $\bar{m}\not\in I^{< B}$ for some numerals $\bar{n},\bar{m}$.
This means that when the rule $(c)^{B}$ is on the main branch,
$\vphi^{B}$ is a $\Sig^{B}$-formula $\bar{m}\in I^{< B}$.

\subsection{Rewritings on logical rules}\label{subsec:toprule}

In this subsection we consider the cases when the top $\Phi$ is a lower sequent of 
a rule $J_{1}$, which is
one of basic rules $(\lor)_{1},(b\exi)^{\Natural}_{1},(\exi)^{\Natural}_{1}, (b\exi)^{\mathcal{O}}_{1}, (\exi)^{\mathcal{O}}_{1}, (I)$ introducing an $\exi$-formula.
Let $J$ denote a $(cut)$ at which the descendant of the principal formula of the rule $J_{1}$ vanishes.
\\

\noindent
{\bf Case 1}. 
Between the top $\Phi$ and  $J$, either there is an $(h)$, or there is a $(cut)$ with its height$=0$.
 Let ${\tt s}:J_{0}$ be the uppermost such rule:
\[
P=
\left.
\begin{array}{c}
\infer*{}
{
 \infer[{\tt s}:J_{0}]{{\tt t}:\Gam; d(H+1)}
 {
  \infer*{\cdots;H+1}
  {
   \infer[J_{1}]{\Phi;\, H_{1}+1}
   {
    \infer*{\cdots;\, H_{1}}{}
    }
  }
 }
} 
\end{array}
\right.
\]
where $H=K+H_{1}$, $d(H+1)=D_{{\tt s}}(H+1)$ if ${\tt s}:J_{0}$ is a $(D)$, 
$d(H+1)=\ome(H+1)$ if ${\tt s}:J_{0}$ is an $(h)$, and $d(H+1)=\vphi(\rk(\vphi);H+1)$ for 
the bounded cut formula $\vphi$ of ${\tt s}:J_{0}$.


For  the minor formula $\tht$ of $J_{1}$, let
\[
P':=
\left.
\begin{array}{c}
\infer*{}
{
 \infer[(p)_{d(H)}]{{\tt t}:\Gam; d(H)\cdot 2}
 {
  \infer[J_{0}]{\tht',\Gam; d(H)}
  {
   \infer*{\cdots;H}
   {
    \infer[(pad)_{0}]{\tht,\Phi; \, H_{1}}
    {
     \infer*{\tht,\Phi;\, H_{1}}{}
     }
    }
  }
 }
} 
\end{array}
\right.
\]
where $\tht'$ is a descendant of $\tht$, which may differ from $\tht$ due to rules $(c)$.
We have $(d(H+1),\emptyset)\to_{\ell}(d(H)\cdot 2,\emptyset)$ by Definition \ref{df:move}.\ref{move:d}.

It is easy to see that $P^{\prime}$ is regular
for the same ordinal assignment to rules $(D)$.


For example, when $(p)_{1}=(b\exi)^{\mathcal{O}}_{1}$ 
with a principal formula $\bar{n}\in I^{<C}:\equiv(\exi a<C(\bar{n}\in I^{a}))\, (C\leq\mu)$
and a minor formula $\bar{n}\in I^{A}$ with $A<C$,
we have $A< B$ in the following figures
by {\bf (p2)} in Definition \ref{df:proofha}.
\[
P=
\left.
\begin{array}{c}
\infer*{}
{
   \infer[J_{0}]{\bar{n}\in I^{< B},\Del_{2};d(H+1)}
   { 
    \infer*{\cdots;H+1}
    {
      \infer[{\tt u}:(c)_{B}^{C}]{\bar{n}\in I^{< B},\Del_{1};H_{0}+1}
     {
       \infer*{\bar{n}\in I^{<C},\Del_{1};H_{0}+1}
       {
         \infer[(b\exi)^{\mathcal{O}}_{1}]{\bar{n}\in I^{<C},\Del_{0};\, H_{1}+1}
         {
          \infer*{\bar{n}\in I^{A},\bar{n}\in I^{<C},\Del_{0};\, H_{1}}{}
          }
        }
       }
      }
    }
} 
\end{array}
\right.
P':=
\left.
\begin{array}{c}
\infer*{}
{
   \infer[(b\exi)^{\mathcal{O}}_{d(H)}]{\bar{n}\in I^{< B},\Del_{2}; d(H)\cdot 2}
   {
    \infer[J_{0}]{\bar{n}\in I^{A},\bar{n}\in I^{<B},\Del_{2};d(H)}
    { 
     \infer*{\cdots;H}
     {
       \infer[(c)_{B}^{C}]{\bar{n}\in I^{A},\bar{n}\in I^{< B},\Del_{1};H_{0}}
      {
        \infer*{\bar{n}\in I^{A},\bar{n}\in I^{<C},\Del_{1};H_{0}}
        {
         \infer[(pad)_{0}]{\bar{n}\in I^{A},\bar{n}\in I^{<C},\Del_{0};\, H_{1}}
         {
          \infer*{\bar{n}\in I^{A},\bar{n}\in I^{<C},\Del_{0};\, H_{1}}{}
          }
          }
        }
       }
      }
    }
} 
\end{array}
\right.
\]
{\bf Case 2}. The cut formula of $J$ is unbounded. 

Due to {\bf  Case 1} in this subsection, there is no $(h)$ nor $(D)$ between $\Phi$ and $J$.
Let $\tht_{1}$ be the unbounded cut formula.
The principal formula of $J_{1}$ is the formula $\tht_{1}$.
Let $\tht_{0}$ be its minor formula.
\[
P=
\left.
\begin{array}{c}  
 \infer*{}
{
 \infer[{\tt s}:J_{0}]{{\tt t}:\Gam_{4};d(H_{4}+1)}
 {
  \infer*{\Gam_{4};H_{4}+1}
  {
   \infer[J]{\Gam,\Del;H_{3}+1}
   {
    \deduce{\Gam,\lnot\tht_{1}; H_{1}}{P_{0}}
   &
    \infer*{\tht_{1},\Del; H_{2}+1}
    {
      \infer[J_{1}]{\tht_{1},\Del_{0};\, H_{0}+1}{\tht_{0}, \tht_{1},\Del_{0};\, H_{0}}
   }
   }
 }
 }
}
\end{array}
\right.
\]
where $H_{2}=H_{5}+H_{0}$ for some $H_{5}$,
$H_{3}=H_{1}+H_{2}$, and ${\tt s}:J_{0}$ is the uppermost rule $(h)$ or $(D)$ below $J$.
$d=\ome$ if $h({\tt t})>0$.
Otherwise $d(H)=D_{{\tt s}}(H)$.
Let $\rk(\tht_{1})=\mu+n_{1}$ with $n_{1}\leq h({\tt t})+1$ by {\bf (p0)}, and $\rk(\tht_{0})=A+n_{0}$.
We have $n_{0}<n_{1}$ if $A=\mu$, and $n_{0}<\ell$.
We have $(\ome(H_{4}+1),\emptyset)\to_{\ell}(\ome(H_{4})\cdot 2,\emptyset)$ by Definition \ref{df:move}.\ref{move:d}, and
by Definition \ref{df:move}.\ref{move:D2}
$(D_{{\tt s}}(H_{4}+1),lb)\to (\vphi_{A+n_{0}+1}(D_{{\tt s}}(H_{4})\cdot 2),lb)$ for
$A\in lb(\tht_{0})\subset L(P)\subset lb$ and $L(P\uarw{\tt s})\ni A\ull_{\mu}\Ome({\tt s})$ by {\bf (p3)}.

Assuming that $\tht_{0}$ is a $\fal$-formula,
let $P'$ be the following with $e(H)\equiv H$ if $d=\ome$, and $e(H)=\vphi_{A+n_{0}}(H)$
if $d=D_{{\tt s}}$ and $\rk(\tht_{0})=A+n_{0}$.
\[
P'=
\left.
\begin{array}{c}  
 \infer*{}
{
\infer{{\tt t}:\Gam_{4};e(d(H_{4})\cdot 2)}
{
 \infer[{\tt s}_{0}]{{\tt t}_{0}:\Gam_{4},\tht_{0};d(H_{4})}
 {
  \infer*{\Gam_{4},\tht_{0};H_{4}}
  {
   \infer{\Gam,\Del,\tht_{0};H_{3}}
   {
    \deduce{\Gam,\lnot\tht_{1}; H_{1}}{P_{0}}
   &
    \infer*{\tht_{1},\tht_{0},\Del; H_{2}}
    {
      \infer{\tht_{1},\tht_{0},\Del_{0};\, H_{0}}{\tht_{0}, \tht_{1},\Del_{0};\, H_{0}}
   }
   }
 }
 }
 &
  \infer[{\tt s}_{1}]{{\tt t}_{1}:\Gam_{4},\lnot\tht_{0};d(H_{4})}
 {
  \infer*{\Gam_{4},\lnot\tht_{0};H_{4}}
  {
   \infer[(pad)_{H_{2}}]{\Gam,\Del,\lnot\tht_{0};H_{3}}
   {
    \deduce{\Gam,\lnot\tht_{0}; H_{1}}{P'_{0}}
   }
 }
 }
}
}
\end{array}
\right.
\]
where $P_{0}^{\prime}$ is obtained from $P_{0}$ by inversion.
When $d=D_{{\tt s}}$, let $\Ome^{\prime}({\tt s}_{i})=\Ome({\tt s})$ for $i=0,1$.
We see that $P^{\prime}$ is regular as follows.
Let $\tht_{1}\equiv\exi a \vphi(a)$ with a bounded formula $\vphi$, and $\tht_{0}\equiv\vphi(C)$.
${\tt s}:J_{0}$ is a rule $(D)$. 
Consider the condition {\bf (p1)} for nodes ${\tt t}_{i}$, $i=0,1$.
Assume $\max\{lh({\tt t}),lq(\tht_{0})\}=lh({\tt t}_{0};P^{\prime})\leq C$.
Then $lh({\tt t})\leq C\in L(P\uarw{\tt t})$. 
We obtain $C\in\Ome_{f}({\tt t})=\Ome({\tt s})=\Ome^{\prime}({\tt s}_{i})=\Ome_{f}^{\prime}({\tt s}_{i})$ by {\bf (p1)}.
\\

\noindent
In what follows assume that the cut formula of $J$ is bounded.

Due to {\bf  Case 1} in this subsection, between $\Phi$ and $J$
there is no $(h)$, $(D)$ nor $(cut)$ with its height$=0$.
\\
{\bf  Case 3}. There is one of rules rule $(h)$ and $(D)$ below $J$.

Let $\tht_{1}$ be the bounded cut formula.
The principal formula of $J_{1}$ is the formula $\tht_{1}$.
Let $\tht_{0}$ be its minor formula.
\[
P=
\left.
\begin{array}{c}  
 \infer*{}
{
 \infer[{\tt s}:J_{0}]{{\tt t}:\Gam_{4};d(H_{4}+1)}
 {
  \infer*{\Gam_{4};H_{4}+1}
  {
   \infer[J]{\Gam,\Del;H_{3}+1}
   {
    \deduce{\Gam,\lnot\tht_{1}; H_{1}}{P_{0}}
   &
    \infer*{\tht_{1},\Del; H_{2}+1}
    {
      \infer[J_{1}]{\tht_{1},\Del_{0};\, H_{0}+1}{\tht_{0}, \tht_{1},\Del_{0};\, H_{0}}
   }
   }
 }
 }
}
\end{array}
\right.
\]
where $H_{2}=H_{5}+H_{0}$ for some $H_{5}$, 
$H_{3}=H_{1}+H_{2}$, and ${\tt s}:J_{0}$ is the uppermost rule $(h)$ or $(D)$ below $J$.
$d=\ome$ if $h({\tt t})>0$.
Otherwise $d(H)=D_{{\tt s}}(H)$.

Assuming that $\tht_{0}$ is a $\fal$-formula,
let $P'$ be the following.
\[
P'=
\left.
\begin{array}{c}  
 \infer*{}
{
\infer{{\tt t}:\Gam_{4};e(d(H_{4})\cdot 2)}
{
 \infer{\Gam_{4},\tht_{0};d(H_{4})}
 {
  \infer*{\Gam_{4},\tht_{0};H_{4}}
  {
   \infer{\Gam,\Del,\tht_{0};H_{3}}
   {
    \deduce{\Gam,\lnot\tht_{1}; H_{1}}{P_{0}}
   &
    \infer*{\tht_{1},\tht_{0},\Del; H_{2}}
    {
      \infer{\tht_{1},\tht_{0},\Del_{0};\, H_{0}}{\tht_{0}, \tht_{1},\Del_{0};\, H_{0}}
   }
   }
 }
 }
 &
  \infer{\Gam_{4},\lnot\tht_{0};d(H_{4})}
 {
  \infer*{\Gam_{4},\lnot\tht_{0};H_{4}}
  {
   \infer[(pad)_{H_{2}}]{\Gam,\Del,\lnot\tht_{0};H_{3}}
   {
    \deduce{\Gam,\lnot\tht_{0}; H_{1}}{P'_{0}}
   }
 }
 }
}
}
\end{array}
\right.
\]
where $e(H)\equiv H$ when $d=\ome$
with $(\ome(H_{4}+1),\emptyset)\to_{\ell} (\ome(H_{4})\cdot 2,\emptyset)$ by Definition \ref{df:move}.\ref{move:d}.
Otherwise $e(H)\equiv \vphi(\rk(\tht_{0});H)$ with $\rk(\tht_{0})=A+n<\mu$.
In this case $(D_{{\tt s}}(H_{4}+1),lb)\to_{\ell} (\vphi_{A+n+1}(D_{{\tt s}}(H_{4})\cdot 2),lb)$
 by Definition \ref{df:move}.\ref{move:D2}
for $A\leq \Ome({\tt s})$ and $n< \ell$, cf.\,(\ref{eq:cprime}).
$A\ull_{\mu} \Ome({\tt s})$ is seen from $A\in lb(\tht_{0})\cup\{0\}$ and {\bf (p3)}.
Under the same assignment of ordinals to rules $(D),(c)$,
we see that $P^{\prime}$ is regular as in {\bf  Case 2} of this subsection.
\\

\noindent
In what follows assume that there is no rule $(h)$ nor $(D)$ below $J$.
\\
{\bf  Case 4}.
$J_{1}$ is one of rules $(b\exi)^{\Natural},(b\exi)^{\mathcal{O}}$ introducing bounded quantifier whose
minor formula $\tht_{0}$ contains a false immediate subformula.
This means for example, $\tht_{1}\equiv(\exi a<B\tht(a))$, $\tht_{0}\equiv(A<B\land\tht(A))$ with $A\not<B$.
Then replace the rule $J_{1}$ by $(pad)_{0}$ by Necrosis,
\[
P=
\left.
\begin{array}{c}
    \infer*{}
    {
      \infer[J_{1}]{\tht_{1},\Del_{0};\, H_{0}+1}
      {
      \infer*[P_{0}]{\tht_{0}, \tht_{1},\Del_{0};\, H_{0}}{}
      }
   }
\end{array}
\right.
\msfiv
P^{\prime}=
\left.
\begin{array}{c}
    \infer*{}
    {
      \infer[(pad)_{0}]{\tht_{1},\Del_{0};\, H_{0}}
      {
      \infer*[P_{0}^{\prime}]{\tht_{1},\Del_{0};\, H_{0}}{}
      }
   }
\end{array}
\right.
\]
$P_{0}^{\prime}$ is obtained from $P_{0}$ by inversion, and eliminating the false prime formula $A<B$.
\\

\noindent
{\bf Case 5}.
Let $\tht_{1}^{\prime}$ be the bounded cut formula.
The principal formula $\tht_{1}$ of $J_{1}$ may differ from the formula $\tht_{1}^{\prime}$ due to rules $(c)$.
Let $\tht_{0}$ be its minor formula.
\[
P=
\left.
\begin{array}{c}  
 \infer*{}
{
   \infer[J]{{\tt t}:\Gam,\Del;\vphi(\rk(\tht_{1}^{\prime});H_{3}+1)}
   {
    \deduce{{\tt u}:\Gam,\lnot\tht_{1}^{\prime}; H_{1}}{P_{0}}
   &
    \infer*{\tht_{1}^{\prime},\Del; H_{2}+1}
    {
      \infer[J_{1}]{\tht_{1},\Del_{0};\, H_{0}+1}{\tht_{0}, \tht_{1},\Del_{0};\, H_{0}}
   }
   }
}
\end{array}
\right.
\]
where $H_{2}=H_{5}+H_{0}$ for some $H_{5}$, and
 $h({\tt t})=0$, $H_{3}=H_{1}+H_{2}$.

Assuming that $\tht_{0}$ is a $\fal$-formula, let $P'$ be the following with
$B+m=\rk(\tht_{0}^{\prime})<\rk(\tht_{1}^{\prime})=A+n$ and
$\mbox{\boldmath$C$}:=\{C\in L(P\uarw{\tt u})\cup lb(\tht_{0}): \max\{lh({\tt t}),lq(\tht_{0}^{\prime})\}\leq C<lq(\tht_{1}^{\prime})\}<\rk(\tht_{1}^{\prime})$.
We have
$(\vphi_{A+n}(H_{3}+1),lb)\to_{\ell}
(\vphi_{B+m}(\vphi_{A+n}(H_{3})+\vphi_{\mbox{\scriptsize \boldmath$C$}}(H_{3})),lb)$ by 
Definition \ref{df:move}.\ref{move:vphi2},
where $B\in lb$, $m< \ell$ by (\ref{eq:cprime}), and either $B\leq A$, or $B=A$ and $m<n$.
\[
P'=
\left.
\begin{array}{c}  
 \infer*{}
{
 \infer{{\tt t}:\Gam,\Del;\vphi(\rk(\tht_{0}^{\prime});\vphi(\rk(\tht_{1}^{\prime});H_{3})+\vphi(\mbox{\boldmath$C$};H_{3})))}
 {
   \infer{\Gam,\Del,\tht_{0}^{\prime};\vphi(\rk(\tht_{1}^{\prime});H_{3})}
   {
    \deduce{\Gam,\lnot\tht_{1}^{\prime}; H_{1}}{P_{0}}
   &
    \infer*{\tht_{1}^{\prime},\Del,\tht_{0}^{\prime}; H_{2}}
    {
      \infer{\tht_{1},\Del_{0},\tht_{0};\, H_{0}}{\tht_{0}, \tht_{1},\Del_{0};\, H_{0}}
   }
   }
  &
  \infer[(rank)_{\mbox{\scriptsize \boldmath$C$}}]{{\tt u}^{\prime}:\Gam,\lnot\tht_{0}^{\prime};\vphi(\mbox{\boldmath$C$};H_{3})}
  {
    \infer[(pad)_{H_{2}}]{\Gam,\lnot\tht_{0}^{\prime};H_{3}}
    {\deduce{\Gam,\lnot\tht_{0}^{\prime}; H_{1}}{P_{0}'}
    }
  }
  }
}
\end{array}
\right.
\]
where $P_{0}'$ is obtained from inversion.
We see that $P^{\prime}$ is regular as follows.
Let us examine the case when $\tht_{1}^{\prime}$ is a formula $\exi a< A\, \tht(a)$, where
$\tht_{1}^{\prime}\equiv\tht_{1}$ unless $\tht_{1}^{\prime}\equiv(n\in I^{<A})$, $\tht_{1}\equiv(n\in I^{<B})$
and there exists a rule $(c)^{B}_{A}$ between $J_{1}$ and $J$.
In the latter case with $B=\mu$, $J_{1}$ is a rule $(\exi)^{\mathcal{O}}$.
Otherwise $J_{1}$ is a rule $(b\exi)^{\mathcal{O}}$.
Assuming $\tht_{1}^{\prime}\not\equiv\tht_{1}$, the witnessing constant
$A_{0}$ for $n\in I^{<B}$, i.e., $\tht_{0}\equiv(A_{0}<B\land n\in I^{A_{0}})$
is smaller than $A$ if $A_{0}<B$, where we can assume that $A_{0}<B$ due to {\bf Case 4}.

First consider the condition {\bf (p1)} for the node ${\tt u}^{\prime}$ in $P^{\prime}$.
Assume $lh({\tt u}^{\prime})=\max\{lh({\tt t}),lq(\tht_{0}^{\prime})\}>lh({\tt t})$.
Let $C\in L(P^{\prime}\uarw{\tt u}^{\prime})$ be a label such that $C\geq lh({\tt u}^{\prime})$.
If $C\geq lq(\tht_{1}^{\prime})\geq A$, then we see $C\in L(P\uarw{\tt u})$
since if $C\in lb(\tht_{0})$ and $C\not\in lb(\tht_{1}^{\prime})$, then $C<lq(\tht_{1}^{\prime})$.

Therefore $C\in\Ome_{f}({\tt u})\subset\Ome^{\prime}_{f}({\tt u}^{\prime})$ by {\bf (p1)} for $P$.
Let $\max\{lh({\tt t}),lq(\tht_{0}^{\prime})\}=lh({\tt u}^{\prime})\leq C<lq(\tht_{1}^{\prime})$. 
Then $C\in \mbox{\boldmath$C$}$, and hence $C\in\Ome^{\prime}_{f}({\tt u}^{\prime})$.

Next the condition {\bf (p2)} for rules $(c)^{C}$, 
and {\bf (p3)} for rules $(D)$ in $P_{0}^{\prime}$.
Let ${\tt v}^{\prime}:J^{\prime}$ be one of such rule, and
$\tht_{1}^{\prime}\equiv(\exi a< A\, \tht(a))$, and $\tht_{0}\equiv(A_{0}<A^{\prime}\land \tht(A_{0}))$,
where either $A^{\prime}=A$ or $A=d_{A^{\prime}}(K)$ for a $K$ due to a rule $(c)^{A^{\prime}}_{A}$.
The constant $A_{0}$ may occur in $P_{0}^{\prime}$ when  the variable $a$ occurs in $\tht(a)$,
although it need to do so in $P_{0}$.
We have $A_{0}<A\leq lq(\tht_{1}^{\prime})$ if $A_{0}<A$.

Consider first the case when $\lnot\tht_{1}^{\prime}\equiv(\fal a<A\lnot\tht(a))$ is in the upper sequent
of the corresponding rule ${\tt v}:J$ in $P$.
Then $L(P\uarw{\tt v})\ni A\ll_{C}K$ if $J$ is a rule $(c)^{C}_{d_{C}(K)}$, and
$A\ull_{\mu}\Ome({\tt v})$ if $J$ is a $(D)$.
In the former case we have $A<C$, and $A_{0}\ll_{C}A$, and in the latter
$A_{0}\ll_{\mu}A<\mu$.
Finally consider the case when $\lnot\tht_{1}^{\prime}\equiv(\fal a<A\lnot\tht(a))$ is in the lower sequent
of the corresponding rule ${\tt v}:(c)^{B}_{A}$ in $P$, and $\fal a<B\lnot\tht(a)$ is in its upper sequent.
Let $A=d_{B}(K)$. We need to show $A_{0}\ll_{B}K$.
We have $A_{0}\ll_{B^{+}}B$ for $L(P\uarw{\tt v})\ni B\ll_{B}K$.
On the other hand we have $A_{0}<d_{B}(K)=A$.
Hence $A_{0}\ll_{B}K$.

\subsection{Rewritings on induction and reflection}\label{subsec:indrfl}
In this subsection we consider the cases when the top $\Phi$ is a lower sequent of 
one of rules $(VJ), (TJ),(Cl)$.
\\
{\bf Case 1}. $\Phi$ is the lower sequent of a $(VJ)$.
\[
P=
\left.
\begin{array}{c}    
\infer*{}
{
\infer[(VJ)]{{\tt t}:\Phi; (H_{2}+1)\cdot mj(n)}
 {
   \infer*{\Phi,\vphi(0); H_{1}}{}
   &
   \infer*{\lnot \vphi(x),\Phi, \vphi(x') ; H_{2}}{}
   &
   \infer*{\lnot\vphi(n),\Phi; H_{3}}{}
  }
 }
\end{array}
\right.
\]
where $H_{2}=H_{1}+H_{3}<\ome$, $\vphi$ is a $\fal$-formula, and $h({\tt t})>0$ by {\bf (p0)}.
\\
{\bf Case 1.1}.
$mj(n)=*_{\ome}$.
Let
\[
P'=
\left.
\begin{array}{c}    
\infer*{}
{
\infer[(VJ)]{\Phi; (H_{2}+1)\cdot(1+n)}
 {
   \infer*{\Phi,\vphi(0); H_{1}}{}
   &
   \infer*{\lnot \vphi(x),\Phi, \vphi(x') ; H_{2}}{}
   &
   \infer*{\lnot\vphi(n),\Phi; H_{3}}{}
  }
 }
\end{array}
\right.
\]
where $\ell> n$ with $P=P[\ell]$ by (\ref{eq:cprime}).
$((H_{2}+1)\cdot *_{\ome},\emptyset)\to_{\ell}((H_{2}+1)\cdot(1+n),\emptyset)$ 
by Definition \ref{df:move}.\ref{move:cdot*}.
\\
{\bf Case 1.2}.
$mj(n)=1+n$.

When $n>0$, let $P'$ be the following with $mj(n-1)=n$. 
$((H_{2}+1)\cdot(1+n),\emptyset)\to_{\ell}((H_{2}+1)\cdot n+H_{2},\emptyset)$
by Definition \ref{df:move}.\ref{move:cdot*}.
 {\small
\[
\infer*{}
{
\infer[(cut)]{\Phi;(H_{2}+1)\cdot n+H_{2}}
{
\infer[(VJ)]{\Phi,\vphi(n); (H_{2}+1)\cdot n}
 {
   \infer*{\Phi,\vphi(0); H_{1}}{}
   &
   \infer*{\lnot \vphi(x),\Phi, \vphi(x') ; H_{2}}{}
   &
   \infer*{\lnot \vphi(n-1),\Phi, \vphi(n) ; H_{2}}{}
  }
&
 \infer[(pad)_{H_{1}}]{\lnot\vphi(n),\Phi;H_{2}}
 {
   \infer*{\lnot\vphi(n),\Phi; H_{3}}{}
  }
}
}
\]
}
If $n=0$, then let $P'$ be the following. 
$((H_{2}+1)\cdot 1,\emptyset)\to_{\ell}(H_{2},\emptyset)$ by Definition \ref{df:move}.\ref{move:cdot*}.
\[
\infer*{}
{
\infer[(cut)]{\Phi; H_{2}}
 {
   \infer*{\Phi,\vphi(0); H_{1}}{}
   &
   \infer*{\lnot\vphi(0),\Phi; H_{3}}{}
  }
 }
 \]
\\
{\bf Case 2}. $\Phi$ is the lower sequent of a $(TJ)$.
\[
P=
\left.
\begin{array}{c}    
\infer*{}
{
 \infer[(TJ)]{{\tt t}:A\not< B,\Phi; (H+H_{\vphi})\cdot mj(B)}
 {
  \infer*{\Phi,\vphi(a),\lnot\fal b< a\, \vphi(b);H}{}
  &
    \infer*{\lnot \vphi(A),\Phi;H_{\vphi}}{}
 }
}
\end{array}
\right.
\]
where $A,B$ are constants and $mj(B)\in\{B,*_{\mu}\}$, $\vphi$ is a $\fal$-formula,
and $h({\tt t})>0$ by {\bf (p0)}.
\\
{\bf Case 2.1}. $A\not< B$: Then $A\not< B$ is a true s.p.f., and $A\not< B,\Gam$ is a stage prime axiom.
This case is reduced to the {\bf Case 2} in subsection \ref{subsec:toppad} by (Necrosis).
\\
{\bf Case 2.2}. $A< B$: Then $A< o(mj(B))$ with $o(B)=B$ and $o(*_{\mu})=\mu$.
\\
We have $A\in lb:=lb[\ell]$ with $P=P[\ell]$.

Let with $A=mj(A)$. $((H+H_{\vphi})\cdot mj(B),lb)\to_{\ell}((H+H_{\vphi})\cdot A+H+H_{\vphi},lb)$ by
Definition \ref{df:move}.\ref{move:cdot*}.
By {\bf (p2)} we have $A\ll_{C}H_{0}$ for any rule ${\tt s}:(c)^{C}_{d_{C}(H_{0})}$ occurring below the $(TJ)$,
${\tt s}\subset_{e}{\tt t}$.
Hence $\Ome({\tt s}_{0};P^{\prime})\ll_{C}H_{0}$ holds for the upper sequent ${\tt s}_{0}$ of ${\tt s}$.
{\bf (p2)} is enjoyed, and $P^{\prime}$ is regular.
{\small
\[
\infer*{}
{
\infer[(cut)]{A\not< B,\Phi; (H+H_{\vphi})\cdot A+H+H_{\vphi}}
{
\infer{\Phi,\fal a< A\, \vphi(a)}
 {
  \infer[(TJ)]{a\not< A,\Phi,\vphi(a);(H+H_{\vphi})\cdot A}
  {
   \infer*{\Phi,\vphi(a),\lnot\fal b< a\, \vphi(b);H}{}
   &
      \infer*{\vphi(a),\lnot\vphi(a),\Phi; H_{\vphi}}{}
   }
  }
&
\hspace{-20mm}
 \infer*[a:=A]{\Phi,\vphi(A),\lnot\fal b< A\, \vphi(b);H}{}
 &
 \infer*{\lnot \vphi(A),\Phi;H_{\vphi}}{}
 }
}
\]
}
{\bf Case 3}. The top ${\tt t}: \Phi$ is the lower sequent of a $(Cl.B)$.

Let $P$ be the following with $\tht(a):\equiv\left(\fal y \calb_{i}(I^{< a},\bar{n},y)\right)$:
\[
 \infer*{}
{
   \infer[(Cl.B)]{{\tt t}: \Phi;\{B^{*}\}(H_{0}+H_{1})}
   {
    \infer*{\Phi,\tht(B);H_{0}}{}
    &
     \infer*[P_{0}]{\lnot R_{B}(a),a\not< B,  \lnot \tht(a),\Phi;H_{1}}{}
     }
}
\]
where $B$ denotes a constant $B\leq\mu$, and $i=0$ when $B\neq\mu$, $B^{*}=\mu$ when $B=\mu$,
and $B^{*}\in\{*_{\mu},B\}$ otherwise.
$\lnot R(B)$ is in $\Phi$.
\\
{\bf Case 3.1}.
$B^{*}=*_{\mu}$:
Then $h({\tt t})>0$ by {\bf (p0)}.
Let $P^{\prime}$ be the following. We have $(\{*_{\mu}\}(H_{0}+H_{1}),lb)\to_{\ell}(\{B\}(H_{0}+H_{1}),lb)$
 by Definition \ref{df:move}.\ref{move:brace1} for $B\in lb[\ell]$.
\[
 \infer*{}
{
   \infer[(Cl.B)]{\Phi;\{B\}(H_{0}+H_{1})}
   {
    \infer*{\Phi,\tht(B);H_{0}}{}
    &
     \infer*[P_{0}]{\lnot R_{B}(a),a\not< B,  \lnot \tht(a),\Phi;H_{1}}{}
     }
}
\]
{\bf Case 3.2}.
$B^{*}=B$ and $R(B)$ is false:
$B$ is a constant$<\mu$, and $\Phi$ is an axiom with the true s.p.f. $\lnot R(B)$.
This case is reduced to the {\bf Case 2} in subsection \ref{subsec:toppad} by (Necrosis).
\\

\noindent
In what follows suppose that $B^{*}=B$ and $R(B)$ is true.
\\
{\bf Case 3.3}.
Either $B=\mu$ and there exists a rule $(h)$ below the top, or
$B\neq\mu$ and below the top, there exists one of rules $(h),(D)$ and $(cut)$
with its cut rank $A+n\geq B$ and $n\leq \ell$ by (\ref{eq:cprime}).

Let ${\tt s}:J$ denote the uppermost such rule, and ${\tt t}:\Gam$ its lower sequent.
Let $d=\ome$ when $J$ is a rule $(h)$, $d=D_{{\tt s}}$ when $J$ is a $(D)$, and
$d=\vphi_{A+n}$ with $B\leq A, n\leq\ell$.

\[
 \infer*{}
{
 \infer[{\tt s}:J]{{\tt t}:\Gam; d(K+\{B\}(H))}
 {
  \infer*{}
  {
   \infer[(Cl.B)]{\Phi;\{B\}(H)}
   {
    \infer*{\Phi,\tht(B);H_{0}}{}
    &
     \infer*{\lnot R_{B}(a),a\not< B,  \lnot \tht(a),\Phi;H_{1}}{}
     }
   }
 }
}
\]
where $H=H_{0}+H_{1}$.

We have by Definition \ref{df:move}.\ref{df:move.brace2} that
$(d(K+\{B\}(H)),lb)\to_{\ell} (\{B\}(d(K+H)\cdot 2),lb)$.

Let $P^{\prime}$ be the following, where $\Ome^{\prime}({\tt s}_{i})=\Ome({\tt s})\, (i=0,1)$
when $d=D_{{\tt s}}$:
\[
 \infer*{}
{
\infer[(Cl.B)]{{\tt t}:\Gam;\{B\}(d(K+H)\cdot 2)}
{
 \infer[{\tt s}_{0}]{\Gam,\tht(B); d^{\prime}(K+H)}
 {
  \infer*{}
  {
   \infer[(pad)_{H_{1}}]{\Phi,\tht(B);H}
   {
    \infer*{\Phi,\tht(B);H_{0}}{}
   }
   }
  }
 &
 \infer[{\tt s}_{1}]{\lnot R_{B}(a),a\not< B,  \lnot \tht(a),\Gam; d^{\prime}(K+H)}
 {
  \infer*{}
  {
    \infer[{}_{H_{0}}(pad)]{\lnot R_{B}(a),a\not< B,  \lnot \tht(a),\Phi;H}
    {
     \infer*{\lnot R_{B}(a),a\not< B,  \lnot \tht(a),\Phi;H_{1}}{}
     }
   }
  }
 }
}
\]
{\bf Case 3.4}. $B=\mu$ and $h({\tt t})=1$ for the top ${\tt t}:\Phi$.
\[
 \infer*{}
{
 \infer[{\tt s}_{D}:(D)]{{\tt s}:\Gam; D(\mbox{\boldmath$C$};K+\{\mu\}(H))}
 {
  \infer*{\Gam;K+\{\mu\}(H)}
  {
  \infer*{{\tt u}:\Gam_{0}; K+\{\mu\}(H)}
  {
   \infer[(Cl.\mu)]{{\tt t}:\Phi;\{\mu\}(H)}
   {
    \infer*{\Phi,\tht(\mu);H_{0}}{}
    &
     \infer*[P_{0}]{\lnot R(a), \lnot \tht(a),\Phi;H_{1}}{}
     }
   }
   }
 }
}
\]
where $H=H_{0}+H_{1}$ and
${\tt u}:\Gam_{0}$ denotes the upper sequent of the uppermost rule $(c)^{\mu}$ below the top {\tt t}
if such a rule exists.
Otherwise ${\tt u}:\Gam_{0}(=\Gam)$ is the upper sequent of the rule $(D)$ below the top {\tt t}.
We have $\Gam_{0}\subset\Sig^{\mu}$ due to subsection \ref{subsec:topbnd}.
Also $\mbox{\boldmath$C$}\subset lb$ such that $L(P\uarw{\tt u})\ull_{\mu} \mbox{\boldmath$C$}$ 
for $lb=lb[\ell]$ and $P=P[\ell]$ by the condition {\bf (p3)} for the rule ${\tt s}_{D}:(D)$.




Let $P'=P[\ell+1]$ be the following:
{\small
\[
 \infer*{}
{
 \infer[(cut)]{{\tt s}:\Gam;\vphi_{A+d_{\cala}}(D(\mbox{\boldmath$C$}\cup\{A\};K+H)\cdot 2)}
 {
  \infer[{\tt s}_{D_{0}}:(D)]{{\tt s}_{0}:\Gam, \tht(A); D(\mbox{\boldmath$C$}\cup\{A\};K+H)}
  {
   \infer*{\Gam,\tht(A);K+H}
   {
    \infer[{\tt u}_{00}:(c)^{\mu}_{A}]{\Gam_{0},\tht(A); K+H}
    {
     \infer*{{\tt u}_{0}:\Gam_{0},\tht(\mu);K+H}
    {
     \infer[(pad)]{\Phi,\tht(\mu);H}
     {
      \infer*{\Phi,\tht(\mu);H_{0}}{}
     }
     }
     }
   }
  }
&
\hspace{-2mm}
 \infer[{\tt s}_{D_{1}}:(D)]{{\tt s}_{1}:\lnot \tht(A),\Gam; D(\mbox{\boldmath$C$}\cup\{A\};K+H)}
 {
 \infer*{\lnot\tht(A),\Gam;K+H}
 {
    \infer*{{\tt u}_{1}:\lnot\tht(A),\Gam_{0};K+H}
   {
    \infer[(pad)]{\lnot \tht(A),\Phi; H}
    {
      \infer*[P_{0}^{\prime}]{ \lnot \tht(A),\Phi;H_{1}}{}
     }
    }
   }
  }
 }
}
\]
}
where $P_{0}^{\prime}$ is obtained from the subproof $P_{0}$ of $P$
by substituting the constant $A$ for the eigenvariable $a$, and eliminating the false formula
$\lnot R(A)$.
The label heights $lh^{\prime}$ for $P^{\prime}$ are defined by 
$lh^{\prime}({\tt s})=lh({\tt s})$ and $lh^{\prime}({\tt s}_{i})=\max\{lh({\tt s}), lq(\tht(A))\}$.
The condition {\bf (p1)} is enjoyed for $P^{\prime}$ since $A\in(D(\mbox{\boldmath$C$}\cup\{A\};K+H))_{f}$.
Let $A=d_{\mu}(\{A_{p}\}(K+H))\in L$ for $A_{p}=\max(\{0\}\cup L(P\uarw{\tt u}))$.
We have 
$\rk(\tht(A))=\rk(\fal y\fal b< A\, \calb_{i}(I^{b},I^{< A},n,y))=A+d_{\cala}$ with
$d_{\cala}\leq\ell$ by (\ref{eq:cprime}).
For the rules $(D)$ in $P^{\prime}$, let
$\Ome^{\prime}({\tt s}_{D_{i}})=\mbox{\boldmath$C$}\cup\{A\}$.
Then we have
$(D(\mbox{\boldmath$C$};K+\{\mu\}(H)),lb)\to_{\ell} (\vphi_{A+d_{\cala}}(D(\mbox{\boldmath$C$}\cup\{A\};K+H)\cdot 2),lb\cup\{A\})$
by (Production) in Definition \ref{df:move}.\ref{move:bracemu}.
The condition {\bf (p3)} is fulfilled for rules ${\tt s}_{D_{i}}:(D)\,(i=0,1)$.

Let 
$lb[\ell+1]=lb^{\prime}=\{A\}\cup lb$.
Consider the condition {\bf (p2)} for the new rule ${\tt u}_{00}:(c)^{\mu}_{A}$.
We have $L(P^{\prime}\uarw{\tt u}_{0})\subset L(P\uarw{\tt u})\ull_{\mu}A_{p}\ll_{\mu}\{A_{p}\}(K+H)$
for $A_{p}<\mu$.

Next let $\infer[{\tt v}_{0}:(c)^{\mu}_{C}]{\Gam_{1}}{{\tt v}:\Gam_{2}}$ with $C=d_{\mu}(H_{0})$
be a rule above {\tt s} in $P$.
Consider the condition {\bf (p2)} for the corresponding rule 
${\tt v}_{0}^{\prime}:(c)^{\mu}_{C}$ in $P^{\prime}$.
{\bf (p2)}  is fulfilled since $A\ll_{\mu}H_{0}$ from 
$L(P\uarw{\tt v}_{0})\ni A_{p}\ll_{\mu}H_{0}$ and $\{A_{p}\}(K+H)<K+\{\mu\}(H)\leq H_{0}$.


Third consider the condition {\bf (p2)} for rules 
$\infer[{\tt v}_{0}:(c)^{B}_{C_{1}}]{\Gam_{1}}{{\tt v}:\Gam_{2}}$ below {\tt s}.
We have $B<\mu$.
Let $C_{1}=d_{B}(K_{1})$, $K_{2}=\Ome({\tt v})$ and $K_{2}^{\prime}=\Ome^{\prime}({\tt v})$.
We have $K_{2}\ull_{B}K_{1}$ and $L(P\uarw{\tt v}_{0})\ni A_{p}\ll_{B}K_{1}$.
We need to show $A\ll_{B}K_{1}$ and
$K_{2}^{\prime}\ll_{B}K_{2}$.
We see 
$\Ome^{\prime}({\tt s})=\vphi_{A+d_{\cala}}(D(\mbox{\boldmath$C$}\cup\{A\};K+H)\cdot 2)\ll_{\mu} 
D(\mbox{\boldmath$C$};K+\{\mu\}(H))=\Ome({\tt s})$ from
$A\ll_{\mu}D(\mbox{\boldmath$C$};K+\{\mu\}(H))$, which in turn follows from 
$A_{p}\ll_{\mu}D(\mbox{\boldmath$C$};K+\{\mu\}(H))$.
Hence $A,K_{2}^{\prime}\ll_{\mu}K_{2}$.
On the other hand we have $L(P\uarw{\tt v}_{0})\ni A_{p}\ll_{B}K_{1}$ and $K_{2}\ull_{B}K_{1}$.
Hence we obtain $A\ll_{B}K_{1}$ and $K_{2}^{\prime}\ll_{B}K_{2}$.

Therefore $(P^{\prime},lb^{\prime})$ is a regular proof.
\\

\noindent
{\bf Case 3.5}. $B\neq\mu$ and below the top, there is no rule $(h),(D)$ and $(cut)$ with its cut rank
$A+n\geq B$.

\[
  \infer*{}
  {
   \infer*{{\tt t}_{1}:\Phi_{1};e(\{B\}(H))}
   {
   \infer[(Cl.B)]{{\tt t}:\Phi;\{B\}(H)}
   {
    \infer*{\Phi,\tht(B);H_{0}}{}
    &
     \infer*[P_{0}]{{\tt t}_{1}:a\not< B,  \lnot \tht(a),\Phi;H_{1}}{}
     }
     }
   }
\]
where $H=H_{0}+H_{1}$, and $B>lh({\tt t})=lh({\tt t}_{1})$ with the lowest ${\tt t}_{1}\subset_{e}{\tt t}$, cf.\,Lemma \ref{lem:rank}.\ref{lem:rank0}.
We have $\Phi_{1}\subset\Sig^{B}$ due to subsection \ref{subsec:topbnd}, and $H\in\calh_{1}$.
Let $A_{p}=\max((L(P\uarw{\tt t}_{1})\cap B)\cup\{0\})$
and
$A=d_{B}(\{A_{p}\}(e(H)))$.
We have $L(P\uarw{\tt t}_{1})\cap B\leq A_{p}\in lb\cup\{0\}$,
and $\rk(\tht(A))=\rk(\fal y\fal b< A\, \calb_{i}(I^{b},I^{< A},n,y))=A+d_{\cala}$ with
$d_{\cala}\leq\ell$ by (\ref{eq:cprime}).

We have
 $(e(\{B\}(H)),lb)\to_{\ell}\left(\vphi_{A+d_{\cala}}(\vphi_{A}(e(H))\cdot 2), lb\cup\{A\}\right)$
by (Production) in Definition \ref{df:move}.\ref{move:braceB}.
Let $P'$ be the following:
\[
 \infer*{}
{
 \infer[(cut)]{{\tt t}:\Phi_{1};\vphi_{A+d_{\cala}}(\vphi_{A}(e(H))\cdot 2)}
 {
 \infer[(rank)_{A}]{{\tt u}_{0}:\Phi_{1}, \tht(A); \vphi_{A}(e(H))}
 {
  \infer[{\tt t}_{0}: (c)^{B}_{A}]{\Phi_{1}, \tht(A); e(H)}
  {
   \infer*{\Phi_{1},\tht(B); e(H)}
   {
     \infer[(pad)]{\Phi,\tht(B);H}
     {
      \infer*{\Phi,\tht(B);H_{0}}{}
     }
     }
     }
  }
&
\hspace{-2mm}
\infer[(rank)_{A}]{{\tt u}_{1}:\lnot \tht(A),\Phi_{1}; \vphi_{A}(e(H))}
{
 \infer{\lnot \tht(A),\Phi_{1}; e(H)}
 {
 \infer*{\lnot\tht(A),\Phi_{1};e(H)}
 {
    \infer[(pad)]{\lnot \tht(A),\Phi; H}
    {
      \infer*[P_{0}^{\prime}]{ \lnot \tht(A),\Phi;H_{1}}{}
     }
   }
   }
  }
 }
}
\]
where $P_{0}^{\prime}$ is obtained from the subproof $P_{0}$ of $P$
by substituting the constant $A$ for the eigenvariable $a$, and eliminating the false formula
$A\not< B$.
The conditions {\bf (p1)} and {\bf (p3)} are enjoyed for $P^{\prime}$.

Consider the condition {\bf (p2)} for the new ${\tt t}_{0}:(c)^{B}_{A}$.
Let $C\in L(P^{\prime}\uarw{\tt t}_{0})$.
We have $C\in L(P^{\prime}\uarw{\tt t}_{0})\subset L(P\uarw{\tt t})$.
$L(P\uarw{\tt t})\cap B\leq A_{p}$ yields $C\ll_{B}\{A_{p}\}(H)$ if $C<B$.
Let $C\geq B>lh({\tt t})=lh({\tt t}_{1})$ for the lowest ${\tt t}_{1}\subset_{e}{\tt t}$.
Then $C\in\Ome_{f}({\tt t}_{1})$ by the condition {\bf (p1)}.
Hence $C\ll_{0}e(H)$, and $C\ll_{B}\{A_{p}\}(e(H))$.

Consider the condition {\bf (p2)} for ${\tt s}:(c)^{C_{0}}_{C_{1}}$ in $P^{\prime}$
other than the new ${\tt t}_{0}:(c)^{B}_{A}$. Let $C_{1}=d_{C_{0}}(K)$.

First let ${\tt t}_{1}\subset_{e}{\tt s}$. Then $C_{0}>B$. We can assume that the formula $a\not<B$
is in the upper sequent ${\tt s}_{0}$ of ${\tt s}$, i.e., $B\in L(P\uarw{\tt s})$.
We obtain $A\ll_{B^{+}}B\ll_{C_{0}}K$, and {\bf (p2)} is enjoyed.

Second let ${\tt s}\subset_{e}{\tt t}_{1}$ and $C_{0}\leq B$.
For $C_{1}=d_{C_{0}}(K)$ it suffices to show $A=d_{B}(\{A_{p}\}(e(H)))\ll_{C_{0}}K$.
We have by {\bf (p2)} that $K^{\prime}\ull_{C_{0}}K$
for $K^{\prime}=\Ome({\tt s}_{0};P)$.
Obviously $\{A_{p}\}(e(H))\ll_{B^{+}}e(\{B\}(H))\ull_{0} K^{\prime}$.
On the other hand we have $B,A_{p},H\ll_{C_{0}}K$ by {\bf (p2)} and $B,A_{p}\in L(P\uarw{\tt t})$.
Hence $A=d_{B}(\{A_{p}\}(e(H)))<d_{B}(K)=C_{1}$, and
$A=d_{B}(\{A_{p}\}(e(H)))\ll_{C_{0}}K$.
\\

This completes a proof of Theorem \ref{th:main}.\ref{th:main3}.
Theorem \ref{th:main}.\ref{th:main1} is proved similarly. Given a proof of a $\Sig_{1}^{0}$-formula
$\exi x \tht(y,x)$, substitute a numeral $\bar{n}$ for the variable $y$, and
add a $(cut)$ with the cut formula $\fal x B(x)\equiv \lnot\exi y \tht(\bar{n},y)$.

We obtain $c=c(P)\geq n+1$ for a proof $P$ of the empty sequent in 
$[\Pi^{0}_{1},\Pi^{0}_{1}]\mbox{{\rm -Fix}} +\fal xB(x)$ .
Let $P_{0}$ be a proof with a $(pad)_{c}$ as the last rule.
Begin to rewrite proofs with $P[c]$ as in this section.
Assuming that $\exi y \tht(\bar{n},y)$ is false, we obtain an infinite path through the tree.

\subsection{Linearity and a theory $[\Pi^{0}_{1},\Pi^{0}_{1}]$-Fixp}\label{sec:linear}

Let $\Phi$-Fixp denote the theory obtained from $\Phi$-Fix by dropping the axiom (\ref{eq:trichotomy}) for
the trichotomy.
Namely in the weakened theory $<$ is supposed to be a \textit{wellfounded} partial order, but the 
\textit{linearity} is not assumed.

It is  easy to see $|a|=|b| \Rarw I^{a}=I^{b}$ without assuming the linearity, 
where $|a|$ denotes the rank $\sup\{|b|+1: b<a\}$.
Therefore $\Phi$-Fixp is supposed to be equivalent to $\Phi$-Fix.
Indeed, the wellfoundedness proofs in \cite{wienpi3} are formalizable in 
$[\Pi^{0}_{1},\Pi^{0}_{1}]$-Fixp. 
Thus theories {\sf KPM}, $[\Pi^{0}_{1},\Pi^{0}_{1}]$-Fix , $[\Pi^{0}_{1},\Pi^{0}_{1}]$-Fixp are 
proof-theoretically equivalent each other
(, i.e., have the same $\Pi^{1}_{1}$-theorems on $\ome$).

Specifically the linearity of the relation $a<b$ was used only in the proof of Theorem 4.4 in \cite{wienpi3}.
Let $\Gam$ denote the operator $\Gam_{2}$ in Definition 4.2 for $Od(\mu)$
in \cite{wienpi3}.
The theorem can be restated as follows:
\blem\label{th:AM} 
$\Gam(I^{a})\ni\alp<\bet\in\Gam(I^{b}) \Rarw \alp\in I^{b}$.
\elem
\bprf 
We show the theorem by induction on the natural sum $|a|\#|b|$ of ordinals (or by main induction on $a$ 
with subsidiary induction on $b$, or vice versa) without assuming the linearity of $<$.

Suppose $\alp<\bet$ and $\alp\in\Gam(I^{a})$, $\bet\in\Gam(I^{b})$.
Then $\alp\in\calg(I^{a})$ by the definition of the operator $\Gam$.
The operator $\calg$ is defined in Definition 3.8.1 in \cite{wienpi3}.

By IH we have
\beqn\label{eq:thA1}
I^{a}|\alp= I^{b}|\alp
\eeqn
Hence by the persistency of $\calg$, cf.\,Lemma 3.9 in \cite{wienpi3}
$\alp\in\calg(I^{a})|\bet=\calg(I^{b})|\bet$.
This suffices to see $\alp\in I^{b}$, cf.\,\cite{wienpi3}.
\eprf
\\

By Lemma \ref{th:AM} we obtain the following equivalence.

\bcor\label{cor:fixp}
The 1-consistency of {\sf KPM} is  equivalent to that of  $[\Pi^{0}_{1},\Pi^{0}_{1}]${\rm -Fixp} over {\sf EA}.
\ecor

We don't need to interpret the relation $A<B$ for labels $A,B\in L$ as in (\ref{eq:less}), 
and there is a chance to replace it by a partial order $A\prec B$, which enjoys
$A\prec B\Rarw o(A)<o(B)$.
Such a partial order $A\prec B$ could be defined through moves $\to_{\ell}$ on hydras
since labels $A,B$ are essentially hydras.
However the trichotomy seems to be indispensable in defining rewritings, e.g., in {\bf Case 3.3}
of subsection \ref{subsec:indrfl}.


\begin{thebibliography}{99}
\bibitem{pntid} T. Arai,  Consistency proof via pointwise induction, 
Arch. Math. Logic 37 (1998), pp. 149-165.


\bibitem{odMahlo} T. Arai, Ordinal diagrams for recursively Mahlo universes, Arch. Math Logic 39 (2000),
pp. 353-391.


\bibitem{ptMahlo} T. Arai, Proof theory for theories of ordinals I:recursively Mahlo ordinals, 
Ann. Pure Appl. Logic  122 (2003), pp. 1-85.


\bibitem{wienpi3} T. Arai, Wellfoundedness proofs by means of non-monotonic inductive definitions I: 
$\Pi^{0}_{2}$-operators, 
Jour. Symb. Logic 69 (2004), pp. 830-850.





\bibitem{Buchholz87} W. Buchholz, An independence result for $(\Pi^1_1-CA)+BI$,
Ann. Pure Appl. Logic 33 (1987), pp. 131-155.

\bibitem{Buchholz90}
W. Buchholz, A note on the ordinal analysis of {\sf KPM},
in Logic Colloquium '90, J. Oikkonen and J. V\"a\"an\"anen (eds.),
Lect. Notes Logic 2, ASL, Cambridge UP, 1993, pp. 1-9.


\bibitem{Buchholz98} W. Buchholz, Relating ordinals to proofs in a perspicuous way, in 
Reflections on the foundations of mathematics (Stanford, CA, 1998),  Lect. Notes Logic 15, ASL, 
Cambridge UP, 2002, pp. 37-59.

\bibitem{KS} L. Kirby and J. Paris, Accessible independence results for Peano arithmetic, Bull. Lon. Math. Soc. 14(1982), pp. 285-293.




\bibitem{Rathjen91}M. Rathjen, Proof-theoretic analysis of $\mbox{KPM}$, Arch. Math. Logic 30 (1991) pp. 377-403.





\bibitem{Richter-Aczel74} W.H. Richter and  P. Aczel, Inductive definitions and reflecting properties of admissible ordinals, Generalized Recursion Theory, Studies in Logic, vol.79, North-Holland, 1974, pp.301-381.


\bibitem{Sk} Th. Skolem, Proof of some theorems on recursively enumerable sets. Notre Dame J. of Formal Logic 3  (1962), pp. 65-74.
\end{thebibliography}
\end{document}